\documentclass[12pt]{paper}
\usepackage{graphicx}
\usepackage{amsmath}
\usepackage{url}
\usepackage{verbatim}
\usepackage{epsf}
\usepackage{latexsym}
\usepackage{array,tabu,multirow,multicol,longtable}
\usepackage{mathrsfs}
\usepackage{authblk}
\usepackage{amsfonts}
\usepackage{amssymb}
\usepackage{moreverb,footmisc,amsthm}%,fnpct}
\usepackage{color}
\usepackage[left=3cm,right=2cm,top=2cm,bottom=2cm]{geometry}
\usepackage{hyperref}
\usepackage[toc,page]{appendix}

\hypersetup{
	colorlinks=true,
	linkcolor=blue,
	urlcolor=magenta,
	citecolor=red,
}
\usepackage{parskip}
\usepackage{algorithm}
\usepackage{algorithmicx}
\usepackage{algpseudocode}
\newcommand*{\QED}{\hfill\ensuremath{\square}}%

\newenvironment{Proof}{\vspace{1.5ex}{\sc Proof}. }{\vspace{2ex} $\QED$}

\newtheorem{Theorem}{Theorem}

\newtheorem{Lemma}{Lemma}
\newtheorem{Proposition}{Proposition}
\newtheorem{Remark}{Remark}

\def\Fc{\mathcal{F}}

\def\Cb{\mathbb{C}}

\def\Rb{\mathbb{R}}

\def\to{\rightarrow}

\def\scal#1#2{
    \left\langle #1,
    #2\right\rangle
    }

\def\mod#1{
    \left|#1
    \right|
    }

\def\norm#1{
    \left\|
    #1\right\|
    }

\def\ud{\,\mathrm{d}}

\begin{document}

\raggedbottom 
\addtolength{\topskip}{0pt plus 10pt}

%-------------------
\title{{\bf \begin{center}
A unified framework for the regularization of final value time-fractional diffusion equation.
\end{center}}}
\author{Walter C. {\sc SIMO TAO LEE}\thanks{Institut de Mathématiques de Toulouse \\
Email: wsimotao@math.univ-toulouse.fr}}
%\date{September 23, 2019}
%-------------------
\maketitle

\abstract 

This paper focuses on the regularization of backward time-fractional diffusion problem on unbounded domain. This problem is well-known to be ill-posed, whence the need of a regularization method in order to recover stable approximate solution. For the problem under consideration, we present a unified framework of regularization which covers some techniques such as Fourier regularization \cite{yang2015fourier}, mollification \cite{van2020mollification} and approximate-inverse \cite{louis1990mollifier}. We investigate a regularization technique with two major advantages: the simplicity of computation of the regularized solution and the avoid of truncation of high frequency components (so as to avoid undesirable oscillation on the resulting approximate-solution). Under classical Sobolev-smoothness conditions, we derive order-optimal error estimates between the approximate solution and the exact solution in the case where both the data and the model are only approximately known. In addition, an order-optimal a-posteriori parameter choice rule based on the Morozov principle is given. Finally, via some numerical experiments in two-dimensional space, we illustrate the efficiency of our regularization approach and we numerically confirm the theoretical convergence rates established in the paper.
  

\section{Introduction}
In this paper, we study the final value problem
\begin{equation}
\label{main equation}
\begin{cases}
\frac{\partial^\gamma u}{\partial t^\gamma}  = \Delta u  & x \in \Rb^n, \,\,\, t \in (0,T)\\
u(x,T) = g(x) & x \in \Rb^n,
\end{cases}
\end{equation}
where we aim at recovering the initial distribution $u(\cdot,0)$ given the final distribution $u(\cdot,T)$. In \eqref{main equation}, $\gamma \in (0,1)$ and $\frac{\partial^\gamma }{\partial t^\gamma}$ denotes the Caputo fractional derivative \cite{kilbas2006theory} defined as 
\begin{equation*}
\label{Caputo Derivative}
\frac{\partial^\gamma u}{\partial t^\gamma} = \frac{1}{\Gamma(1-\gamma)} \int_0^t \frac{u'(s)}{(t-s)^\gamma} \ud s,
\end{equation*}
where $\Gamma(\cdot)$ is nothing but the Gamma function : $\Gamma(z) = \int_0^{+\infty} t^{z-1} e^{-t} \ud t$.

Time-fractional diffusion equations usually model sub-diffusion processes such as slow and anomalous diffusion processes which failed to be described by classical diffusion models \cite{ balakrishnan1985anomalous, chechkin2005fractional, metzler2000random, podlubnv1999fractional}. Due to the high diversity of such phenomena which are not properly modeled by classical diffusion, time-fractional diffusion problems have gained much attention in last decades. Beyond applications in diffusion processes, time-fractional equation \eqref{main equation} has also been applied to image de-blurring \cite{wang2013total} where the time-fractional derivative allows to capture the memory effect in image blurring.

It is well-known that the ill-posedness of equation \eqref{main equation} comes from the irreversibility of time of the diffusion equation, which is caused by the very smoothing property of the forward diffusion. As a result, very small perturbation of the final distribution $u(\cdot,T)$ may cause arbitrary large error in the initial distribution $u(\cdot,0)$. Hence, a regularization method is crucial in order to recover stable approximate of the initial distribution $u(\cdot,0)$. In this regard, many regularization methods have been applied to final value time-fractional diffusion equation. Let us mention the mollification method \cite{van2020mollification, yang2014mollification}, Fourier regularization \cite{xiong2012inverse,yang2015fourier}, the method of quasi-reversibility \cite{liu2010backward}, Tikhonov method \cite{wang2015optimal}, total variation regularization \cite{wang2013total}, boundary condition regularization \cite{yang2013solving}, non-local boundary value method \cite{hao2019stability}, truncation method \cite{wang2012data}. Yet, the set of regularization methods applied to backward time-fractional diffusion equation still presents some sparsity, especially compared to the set regularization methods for backward classical diffusion problems. 

In this paper, we describe how in the context of regularization of the final value time-fractional diffusion equation \eqref{main equation}, the Fourier regularization \cite{yang2015fourier} and mollification \cite{van2020mollification} are nothing but examples of \textit{approximate-inverse} \cite{louis1990mollifier} regularization. Next, we investigate a regularization technique which yields a better trade-off between stability and accuracy compared to the Fourier regularization \cite{yang2015fourier} and the mollification technique of Van Duc N. et al \cite{van2020mollification}. We consider noisy setting where $u(\cdot,T)$ is approximated by a noisy data $g^\delta$ satisfying 
\begin{equation*}
\label{nois data}
|| u(\cdot,T) - g^\delta ||_{L^2(\Rb^n)} \leq \delta,
\end{equation*}
and we derive order-optimal convergence rates between our approximate solution and the exact solution $u(\cdot,0)$ under classical Sobolev smoothness condition 
\begin{equation}
\label{smoothness cond}
u(\cdot,0) \in H^p(\Rb^n) \quad \text{with} \quad ||u(\cdot,0)||_{H^p(\Rb^n)} \leq E, \quad p>0.
\end{equation}
We also provide error estimates under the more realistic setting where both the data $u(\cdot,T)$ and the forward diffusion operator are only approximately known. The motivation here being that, in practice, the Mittag-Leffler function which plays a major role in the resolution of equation \eqref{main equation} can only be approximated in practice.

The outline of this article is as follows: 

In Section \ref{section regularization}, we discuss existence of solution of equation \eqref{main equation} and reformulate the equation into an operator equation of the form $ A u(\cdot,0) = g $ in which $A$ is a bounded linear operator on $L^2(\Rb^n)$. We present key estimates necessary for the regularization analysis and illustrate the ill-posedness of recovering $u(\cdot,0)$ from $g$. Next, we introduce the framework of regularization which includes Fourier regularization, mollification, and approximate inverse. At last, we introduce our regularization approach. 

Section \ref{section error estimate} deals with error estimates and order-optimality of our regularization technique under the smoothness condition \eqref{smoothness cond}. In this section, we derive error estimates between the approximate solution and the exact solution $u(\cdot,0)$ in Sobolev spaces $H^l(\Rb^n)$ with $l \geq 0$. We also give error estimates for the approximation of early distribution $u(\cdot,t)$ with $t \in (0,T)$. We end the section by presenting analogous error estimates for the case where both the data $g$ and the forward diffusion operator $A$ are only approximately known.

Section \ref{section par choice rule} is devoted to parameter selection rules which is a critical step in the application of a regularization method. Here we propose a Morozov-like a-posteriori parameter choice rule leading to order-optimal convergence rates under smoothness condition \eqref{smoothness cond}. We also present analogous error estimates as that obtained in Section \ref{section error estimate} corresponding to the a-posteriori rule prescribed.

Finally, we study four numerical examples in Section \ref{section numerical experiments} to illustrate the effectiveness of the regularization approach coupled with the parameter choice rule described in Section \ref{section par choice rule}. Moreover, in this Section, we also carry out a numerical convergence rates analysis in order to confirm the theoretical convergence rates given in Section \ref{section error estimate}.

In the sequel, $||f||$ or $||f||_{L^2}$ always refers to the $L^2$-norm of the function $f$ on $\Rb^n$, $||f||_{H^p}$ denotes the Sobolev norm of $f$ on $\Rb^n$ and $||| \cdot ||||$ denotes operator norm of a bounded linear mapping. Throughout the paper, $\widehat{f}$ or $\Fc(f)$ (resp. $\Fc^{-1}(f)$) denotes the Fourier (resp. inverse Fourier) transform of the function $f$ defined as
\begin{equation*}
\label{def Fourier transform}
\widehat{f}(\xi) = \Fc(f)(\xi) = \frac{1}{\sqrt{2 \pi}^{n}}
\int_{\Rb^n} f(x)e^{- i x\cdot\xi} \ud x, \quad \Fc^{-1}(f)(x) = \frac{1}{\sqrt{2 \pi}^{n}}
\int_{\Rb^n} f(\xi)e^{ i x\cdot\xi} \ud \xi, \quad \xi, x \in \Rb^n.
\end{equation*}
%%---------------------------------------------------------------
%%---------------------------------------------------------------
\section{Regularization}\label{section regularization}

Let us start by the following result about the existence and uniqueness of solution of equation \eqref{main equation}.
\begin{Proposition}
\label{Prop Existence weak solution}
For all $\gamma \in (0,1)$ and $g \in H^2(\Rb^n)$, Problem \eqref{main equation} admits a unique weak solution $u \in C([0,T],L^2(\Rb^n)) \cap  C((0,T],H^2(\Rb^n))$. That is, the first equation in \eqref{main equation} holds in $L^2(\Rb^n)$ for all $t\in (0,T)$ and $u(\cdot,t) \in H^2(\Rb^n)$ for all $t \in(0,T)$ with
$$
\lim_{t \to T} || u(\cdot,t) - g ||_{H^2} = 0.
$$
\end{Proposition}
Proposition \ref{Prop Existence weak solution} is merely generalization of \cite[Lemma 2.2]{yang2015fourier} where only the case $n=1$ is considered.
The idea of the proof is merely to check that the formal solution defined by \eqref{relation data and solution} is the weak solution.

Now, let us define the framework that we will consider for the regularization of problem \eqref{main equation}. Consider a data $g \in H^2(\Rb^n)$, by applying the Fourier transform in \eqref{main equation} with respect to variable $x$, we get
\begin{equation}
\label{fourier equation}
\begin{cases}
\frac{\partial^\gamma \widehat{u}(\xi,t)}{\partial t^\gamma} = -  |\xi|^{2} \widehat{u}(\xi,t) & \xi \in \Rb^n, \,\,\, t \in (0,T)\\
\widehat{u}(\xi,T) = \widehat{g}(\xi) & \xi \in \Rb^n.
\end{cases}
\end{equation}
By applying the Laplace transform with respect to variable $t$ in \eqref{fourier equation}, one gets
\begin{equation}
\label{link equation data and sol in fourier domain}
\begin{cases}
\widehat{u}(\xi,t) = \widehat{u}(\xi,0) E_{\gamma,1}(-|\xi|^2 t^\gamma) & \xi \in \Rb^n \\
\widehat{u}(\xi,T) = \widehat{g}(\xi) & \xi \in \Rb^n,
\end{cases}
\end{equation}
where $E_{\gamma,1}$ is the Mittag Leffler function \cite{podlubnv1999fractional} defined as
\begin{equation*}
\label{def mittag-Leffler function}
E_{\gamma,1}(z) = \sum_{k=0}^{+\infty} \frac{z^k}{\Gamma(\gamma k  + 1)}, \quad z \in \Cb.
\end{equation*}
From \eqref{link equation data and sol in fourier domain}, we can deduce the following relation between the solution $u(\cdot,0)$ and the data $g$ in the frequency domain:
\begin{equation}
\label{relation data and solution}
\widehat{u}(\xi,0) = \frac{\widehat{g}(\xi)}{E_{\gamma,1}(-|\xi|^2 T^\gamma)}.
\end{equation}
Moreover, we can also derive the following relation for early distribution $u(\cdot,t)$ with $t \in (0,T)$
\begin{equation}
\label{relation data and early distribution}
\forall t \in (0,T), \quad \widehat{u}(\xi,t) = \frac{E_{\gamma,1}(-|\xi|^2 t^\gamma)}{E_{\gamma,1}(-|\xi|^2 T^\gamma)}\widehat{g}(\xi).
\end{equation}
From equations \eqref{relation data and solution} and \eqref{relation data and early distribution}, we can see that the Mittag Leffler function $E_{\gamma,1}$ plays an important role in time-fractional diffusion equation \eqref{main equation}. Hence, let us recall some key estimates about the function $E_{\gamma,1}$ that will be repeatedly used in the sequel.

\begin{Lemma}
\label{Lemma esti Mittag Leffler func}
Let $\gamma \in [\gamma_0,\gamma_1] \subset (0,1)$, there exists constants $C_1$ and $C_2$ depending only on $\gamma_0$ and $\gamma_1$ such that 
\begin{equation}
\label{key estimate E gamma}
\forall x \leq 0, \quad \frac{C_1}{\Gamma(1-\gamma)} \frac{1}{1-x} \leq E_{\gamma,1} (x) \leq \frac{C_2}{\Gamma(1-\gamma)} \frac{1}{1-x}.
\end{equation}
\end{Lemma}
For a proof of Lemma \ref{Lemma esti Mittag Leffler func}, see \cite[Lemma 2.1]{yang2015fourier}.
From Lemma \ref{Lemma esti Mittag Leffler func}, we can easily derive the next Lemma.
\begin{Lemma}
\label{Lemma E gamma}
Assume $\gamma \in [\gamma_0,\gamma_1] \subset (0,1)$, then for every $\xi \in \Rb^n$ and $t \in (0,T]$, 
\begin{equation}
\label{bounds E gamma Xi t }
\frac{1}{(1 \vee t^\gamma)} \left( \frac{C_1}{\Gamma(1-\gamma)} \frac{1}{1+|\xi|^2} \right) \leq E_{\gamma,1} (-|\xi|^2 t^\gamma) \leq \frac{1}{(1 \wedge t^\gamma)} \left(\frac{C_2}{\Gamma(1-\gamma)} \frac{1}{1+|\xi|^2} \right).
\end{equation}
and
\begin{equation}
\label{bounds fract E gamma Xi t }
\frac{C_1}{C_2} \leq \frac{E_{\gamma,1} (-|\xi|^2 t^\gamma)}{E_{\gamma,1} (-|\xi|^2 T^\gamma)} \leq \frac{C_2}{C_1} \left( \frac{T}{t}\right)^\gamma.
\end{equation}
In \eqref{bounds E gamma Xi t }, $\vee$ denotes the maximum while $\wedge$ denotes the minimum, that is, $1 \vee t^\gamma = \max \,\,\{1, t^\gamma\}$ and $1 \wedge t^\gamma = \min\,\, \{1, t^\gamma\}$.
\end{Lemma}

From \eqref{relation data and early distribution} and \eqref{bounds fract E gamma Xi t }, we get that for every $t \in (0,T)$, 
$
||u(\cdot,t)||_{L^2} \leq (C_2/C_1) \left( T/t\right)^\gamma ||g||_{L^2}.
$
which implies that the problem of recovering $u(\cdot,t)$ from $g$ is actually well posed. However, that of recovering $u(\cdot,0)$ is ill-posed. Indeed, from \eqref{relation data and solution} and \eqref{bounds E gamma Xi t }, we get that
\begin{equation}
\label{illus ill-posedness}
||u(\xi,0) || \geq C_\gamma (1 + |\xi|^2) |\widehat{g}(\xi)|, \quad \text{with} \quad C_\gamma = \Gamma(1-\gamma)(1 \wedge T^\gamma)/C_2.
\end{equation}
From \eqref{illus ill-posedness}, we see that very small perturbations in high frequencies in the data $g$ leads to arbitrary large errors in the solution $u(\cdot,0)$. Therefore, one needs a regularization method to recover stable estimates of $u(\cdot,0)$.

\begin{Remark}
From \eqref{relation data and solution} and \eqref{bounds E gamma Xi t }, we can also derive that 
\begin{equation}
\label{middly ill-posedness}
||u(\cdot,t)||_{L^2} \leq  D_\gamma (1 + |\xi|^2) |\widehat{g}(\xi)|, \quad \text{with} \quad D_\gamma = \Gamma(1-\gamma)(1 \vee T^\gamma)/C_1.
\end{equation}
Estimate \eqref{middly ill-posedness} illustrates the fact that the backward time-fractional diffusion equation is less ill-posed (mildly ill-posed) on the contrary to the classical backward diffusion equation which is exponentially ill-posed. This is actually due to the asymptotic slow decay of the Mittag Leffler function $E_{\gamma,1}(-|\xi|^2 )$ compared to $\exp(-|\xi|^2) = E_{1,1}(-|\xi|^2)$.
\end{Remark}
From \eqref{relation data and solution}, we can reformulate equation \eqref{main equation} into an operator equation 
\begin{equation}
\label{operator equation}
A\, u(\cdot,0) = g.
\end{equation}
where $A: L^2(\Rb^n) \to L^2(\Rb^n)$ is the linear forward diffusion operator operator which maps the initial distribution $u(\cdot,0)$ to the final distribution $g$, that is,
\begin{equation}
\label{def operator A pb}
A = \Fc^{-1} \left( E_{\gamma,1}(-|\xi|^2 T^\gamma)  \right) \Fc.
\end{equation}

In the sequel, given a data $g \in L^2(\Rb^n)$, we aim at recovering $u(\cdot,0) \in L^2(\Rb^n)$. Let $\varphi$ be a smooth real-valued function in $L^1(\Rb^n)$ satisfying $\int_{\Rb^n} \varphi(x) \ud x =1$. It is well-known that the family of functions $(\varphi_\alpha)_{\alpha>0}$ defined by 
\begin{equation*}
\label{def phi beta from phi}
\forall x\in \Rb^n,\quad \varphi_\alpha(x) := \frac{1}{\alpha^n}\varphi\left(\frac{x}{\alpha}\right),
\end{equation*}
 satisfies
\begin{equation}
\label{prop phi beta}
\forall f \in L^2(\Rb^n), \quad \varphi_\alpha \star f \to f \quad \text{in} \quad L^2(\Rb^n) \quad \text{as} \quad \alpha \downarrow 0,
\end{equation}
where $\varphi_\alpha \star f$ is nothing but the convolution of the functions $\varphi_\alpha$ and $f$ defined as
$
\left( \varphi_\alpha \star f \right)(x) = \int_{\Rb^n} \varphi_\alpha(x-y)f(y) \ud y.
$
For $\alpha>0$, let $M_\alpha$ be the mollifier operator defined by
\begin{equation}
\label{def oper Cbeta}
\forall \alpha > 0, \,\, \forall f \in L^2(\Rb^n), \quad M_\alpha f = \varphi_\alpha \star f.
\end{equation}
From \eqref{prop phi beta}, we see that the family of operators $(M_\alpha)_{\alpha>0}$ is an approximation of unity  in $L^2(\Rb^n)$, that is,
\begin{equation}
\label{conver M alpha}
\forall f \in L^2(\Rb^n), \quad M_\alpha f \to f \quad \text{in} \quad L^2(\Rb^n) \quad \text{as} \quad \alpha\downarrow 0 .
\end{equation}

Let $u_\alpha$ be the solution of the equation
\begin{equation}
\label{equation u_alpha}
\begin{cases}
\frac{\partial^\gamma u}{\partial t^\gamma}  = \Delta u  & x \in \Rb^n, \,\,\, t \in (0,T)\\
u(x,T) = (M_\alpha g)(x) & x \in \Rb^n.
\end{cases}
\end{equation}
From \eqref{relation data and solution} and \eqref{relation data and early distribution}, by replacing $g$ by $M_\alpha g$, and using the fact that the $\widehat{M_\alpha g}(\xi) = \sqrt{2 \pi}^n \widehat{\varphi_\alpha}( \xi) \widehat{g}(\xi)$, one gets
\begin{equation}
\label{express u alpha in freq domain}
\begin{cases}
\widehat{u_\alpha}(\xi,0) =  \frac{\sqrt{2 \pi}^n}{E_{\gamma,1}(-|\xi|^2 T^\gamma)} \widehat{\varphi_{\alpha}}(\xi) \widehat{g}(\xi) = \sqrt{2 \pi}^n \widehat{\varphi_{\alpha}}(\xi) \widehat{u}(\xi,0) \\
\widehat{u_\alpha}(\xi,t) = \sqrt{2 \pi}^n \frac{E_{\gamma,1}(-|\xi|^2 t^\gamma)}{E_{\gamma,1}(-|\xi|^2 T^\gamma)} \widehat{\varphi_{\alpha}}(\xi) \widehat{g}(\xi) = \sqrt{2 \pi}^n \widehat{\varphi_{\alpha}}(\xi) \widehat{u}(\xi,t), & t \in (0,T).
\end{cases}
\end{equation}
which yields that
\begin{equation}
\label{eq 00}
\forall t \in [0,T] \quad u_\alpha(\cdot,t) = M_\alpha u(\cdot,t).
\end{equation}

\begin{Proposition}
\label{Prop Reg method}
Assume that there exists a function $\varrho: \Rb_+^* \to \Rb_+$ such that the mollifier kernel $\varphi$ verifies 
\begin{equation}
\label{cond reg method}
\forall \alpha>0, \quad \frac{\widehat{\varphi}(\alpha \xi) }{E_{\gamma,1}(-|\xi|^2 T^\gamma)}  \leq \varrho(\alpha),
\end{equation}
Then the family $(u_\alpha)_{\alpha>0}$ solution of equation \eqref{equation u_alpha} defines a regularization method for problem \eqref{main equation} in $L^2(\Rb^n)$.
\end{Proposition}
\begin{Proof} 
From \eqref{express u alpha in freq domain}, noticing that $\widehat{\varphi_\alpha}(\xi) = \widehat{\varphi}(\alpha \xi)$ and using the fact that the function $E_{\gamma,1}$ is increasing on $\Rb_{-}$ with $E_{\gamma,1}(0) = 1$, we can see that \eqref{cond reg method} implies that for every $t \in [0,T]$ and $\alpha >0$, the mapping $R_{\alpha,t}:  L^2(\Rb^n) \to L^2(\Rb^n)$ which maps the data $g$ to $ u_\alpha(\cdot,t)$ is bounded with $|||R_{\alpha,t}||| \leq \sqrt{2 \pi}^n\varrho(\alpha)$. Moreover, from \eqref{conver M alpha} and \eqref{eq 00}, we deduce that for every $t \in [0,T]$, $u_\alpha(\cdot,t)$ converges to $u(\cdot,t)$ in $L^2(\Rb^n)$ as $\alpha$ goes to $0$.
\end{Proof}\\
From Proposition \ref{Prop Reg method}, we can see that, choosing a kernel $\varphi$ which satisfies condition \eqref{cond reg method} allows to defines a regularization method for equation \eqref{main equation}. Now, let us show that the family of regularization methods defined in this way actually coincides with \textit{approximate-inverse} introduced by Louis and Maass \cite{louis1990mollifier}. 

From the first equation in \eqref{express u alpha in freq domain}, we can derive that
\begin{equation}
\label{exp u alpha 0}
u_\alpha(x,0) =  \int_{\Rb^n} e^{i x \xi}  \frac{\widehat{\varphi_{\alpha}}(\xi) \widehat{g}(\xi)}{E_{\gamma,1}(-|\xi|^2 T^\gamma)} \ud \xi .
\end{equation}
From \eqref{eq 00} and \eqref{exp u alpha 0}, we can reformulate $u_\alpha(x,0)$ as follows
\begin{equation}
\label{simil approx inverse}
\begin{cases}
\vspace{0.25cm}
\forall x \in \Rb^n, \quad u_\alpha(x,0) = \scal{e_\alpha(x,\cdot)}{u(\cdot,0)}_{L^2}, & \text{with} \quad e_\alpha(x,y) = \varphi_\alpha(x-y) \\
\forall x \in \Rb^n, \quad u_\alpha(x,0) = \scal{v_{x,\alpha}}{g}_{L^2},  & \text{with} \quad v_{x,\alpha} = \Fc^{-1} \left(  \frac{e^{i x \xi} \,\widehat{\varphi_{\alpha}}(\xi)}{E_{\gamma,1}(-|\xi|^2 T^\gamma)} \right) .
\end{cases}
\end{equation}
 Moreover, one can easily check that $v_{x,\alpha}$ is nothing but the solution of the adjoint equation 
$$ 
A^* f = e_\alpha(x,\cdot).
$$
Hence, we deduce that the solution $u_\alpha$ of equation \eqref{equation u_alpha} actually corresponds to the \textit{approximate-inverse} \cite{louis1990mollifier} regularized solution of equation \eqref{operator equation}, the pair $(v_{x,\alpha},e_\alpha)$ being what is usually called $\textit{reconstruction kernel - mollifier}$. 

This setting of regularization we just described encompasses many regularization methods that appears distinctively in the literature of the regularization of the final value time-fractional diffusion equation. Each regularization method being a particular choice of the mollifier kernel $\varphi$. 

For the Fourier regularization \cite{yang2015fourier} (where $n =1$), we have 
\begin{equation}
\label{app sol Fourier reg}
\widehat{u_{\xi_{max}}}(\xi,t) = \chi_{[-\xi_{max},\xi_{max}]}(\xi) \frac{E_{\gamma,1}(-|\xi|^2 t^\gamma)}{E_{\gamma,1}(-|\xi|^2 T^\gamma)} \widehat{g}(\xi),
\end{equation}
where $\chi_{\Omega}$ denotes the characteristic function of the set $\Omega$  equal to $1$ on $\Omega$ and $0$ elsewhere. By comparing \eqref{app sol Fourier reg} and \eqref{express u alpha in freq domain}, we readily get that
\begin{equation}
\label{ee 11}
\alpha = 1/\xi_{max}, \quad \text{and} \quad \varphi(x) = \Fc^{-1} \left(\frac{1}{ \sqrt{2 \pi}} \chi_{[-1,1]} \right) = \frac{\sin(x)}{\pi x}.
\end{equation}
 In this case, condition \eqref{cond reg method} merely reads 
$$
\frac{\chi_{[-1,1]}(\xi/\xi_{max})}{E_{\gamma,1}(-|\xi|^2 T^\gamma)}  \leq \varrho\left( \frac{1}{\xi_{max}} \right).
$$
Using \eqref{bounds E gamma Xi t }, we can see that this condition is satisfies with 
$
\varrho(\alpha) = C (1 + 1/\alpha^2)$ where $C = (1\vee T^\gamma)\Gamma(1-\gamma) /C_1$. 

For the mollification method of N. Van Duc et al. \cite{van2020mollification}, the mollifier operator $M_\alpha$ is denoted $S_\nu$ where $S_\nu$ is the convolution by the so called \textit{Dirichlet kernel} $D_\nu$ defined as 
$$
D_\nu(x) = \frac{1}{\pi^n} \prod_{j=1}^n \frac{\sin(\nu x_j)}{x_j}, \quad \text{with} \quad \nu >0 \,\, \text{and} \,\, x \in \Rb^n.
$$
Hence we can deduce that, this merely corresponds to
\begin{equation}
\label{ee 22}
\alpha = 1/\nu, \quad \text{and} \quad \varphi(x) =  \prod_{j=1}^n \frac{\sin(x_j)}{ \pi x_j}.
\end{equation}
Given that the Fourier transform of the kernel $D_ \nu$ is given by
\begin{equation}
\label{FT kernel Van duc}
\Fc \left(D_\nu \right)(\xi) = \chi_{\Lambda}(\xi), \quad \text{with} \quad\Lambda = \{x \in \Rb^n \, :\, |x_j| \leq \nu, \, \,\, j = 1,...,n \},
\end{equation} 
we see that condition \eqref{cond reg method} merely reads
$$
\frac{\chi_{\Lambda}(\xi/\nu)}{E_{\gamma,1}(-|\xi|^2 T^\gamma)}  \leq \varrho\left( \frac{1}{\nu} \right),
$$
which is fulfilled with
$\varrho(\alpha) = C (1 + \sqrt{n}/\alpha^2)$ where $C = (1\vee T^\gamma)\Gamma(1-\gamma) /C_1$.

By the way, from \eqref{ee 11} and \eqref{ee 22}, we can see that the Fourier regularization and the mollification approach of N. Van Duc et al. actually coincides, the latter approach being a generalization of the former to $n-$dimensional case.
From \eqref{FT kernel Van duc}, we can conclude that both regularization approaches are nothing but truncation methods. That is, the regularization is done by merely throwing away high frequency components of the data, which are responsible of the ill-posedness, and conserving unchanged the remaining frequency components. In order word, these two methods can be regarded as \textit{spectral cut-off} methods. It is important to notice that though high frequency components are responsible of ill-posedness, nevertheless, the still carry non-negligible information on the sought solution. Therefore, it is desirable to apply a regularization which do not suppress high frequency components but which applies much regularization to those components compared to low frequency components. Let us point out that mere truncation of high frequency components usually entails Gibbs phenomena and oscillation of the approximate solution which should be avoided as far as possible. This is actually possible by choosing a kernel $\varphi$ whose Fourier transform is supported on the whole domain $\Rb^n$.

Now on, let us consider a mollifier kernel $\varphi$ defined by
\begin{equation}
\label{def of our molllif kernel}
\widehat{\varphi}(\xi) = \frac{1}{\sqrt{2 \pi}^n} \exp(- \tau|\xi|^s),\,\,\, \tau >0,\,\,\, s>0, \quad \text{i.e.} \quad \varphi = \frac{1}{\sqrt{2 \pi}^n} \Fc^{-1} \left( \exp(- \tau|\xi|^s) \right),
\end{equation}
where $\tau$ and $s$ are two free positive parameters.
From \eqref{def of our molllif kernel}, we can see that $\varphi \in L^1(\Rb^n)\cap L^2(\Rb^n)$ and satisfies $\int_{\Rb^n} \varphi(x) \ud x = \sqrt{2 \pi}^n \widehat{\varphi}(0) = 1$. 

\begin{Lemma}
\label{Lemma bound our kernel}
Let $b$ and $d$ be two positive numbers, consider the function $f_{b,d}(x) = (1+ x) e^{- b x^d}$. Then there exists a constant $C$ depending only on $d$ such that
\begin{equation}
\label{key bound}
\sup_{x\geq 0} f_{b,d}(x) \leq \frac{C}{b^{1/d}} \quad \text{as} \quad b \downarrow 0.
\end{equation}
\end{Lemma}

The proof of Lemma \ref{Lemma bound our kernel} is deferred to appendix. Lemma \ref{Lemma bound our kernel} will help us to prove that the kernel $\varphi$ given by \eqref{def of our molllif kernel} allows us to define a regularization method for equation \eqref{main equation}.

\begin{Proposition}
\label{Prop our reg method}
Let $M_\alpha$ be the mollifier operator defined by \eqref{def oper Cbeta} with the kernel $\varphi$ given in \eqref{def of our molllif kernel} with $\tau$ and $s$ being two positive numbers. Then the family $(u_\alpha)_{\alpha>0}$ of function $u_\alpha$ solution of equation \eqref{equation u_alpha} defines a regularization method for equation \eqref{main equation}.
\end{Proposition}

\begin{Proof}
In view of Proposition \ref{Prop Reg method}, it suffices to prove that the kernel $\varphi$ given in \eqref{def of our molllif kernel} verifies \eqref{cond reg method}. By considering \eqref{def of our molllif kernel} and estimate \eqref{bounds E gamma Xi t }, we have
\begin{equation}
\label{ee 0011}
\frac{\widehat{\varphi}(\alpha \xi) }{E_{\gamma,1}(-|\xi|^2 T^\gamma)}  \leq C (1 + |\xi|^2) \exp(-\tau \alpha^s |\xi|^s), \quad \text{with} \quad C = \Gamma(1-\gamma)(1 \vee T^\gamma) /\sqrt{2 \pi}^n C_1.
\end{equation}
The right hand side in \eqref{ee 0011} is nothing but $ C f_{b,d}(|\xi|^2)$ with $b = \tau \alpha^{s}$ and $d = s/2$. Hence from \eqref{key bound}, we deduce that there exists a constant $C$ independent on $\alpha$ such that
\begin{equation*}
\forall \xi \in \Rb^n, \quad \frac{\widehat{\varphi}(\alpha \xi) }{E_{\gamma,1}(-|\xi|^2 T^\gamma)}  \leq \frac{C}{\alpha^2} \quad \text{as} \quad \alpha \to 0,
\end{equation*}
whence \eqref{cond reg method} with $\varrho(\alpha)  = C/\alpha^2$.
\end{Proof}\\
\begin{Remark}
By defining the mollifier kernel $\varphi$ as in \eqref{def of our molllif kernel}, we can see that the regularization technique induces a more suitable treatment of frequency components. Indeed, with our choice of mollifier kernel, the amount of regularization smoothly depends on the magnitude of the frequency components: The higher the frequency, the stronger the regularization applied, and similarly, the lower the frequency, the lower the regularization applied. This is actually desirable for a regularization method given that as the frequency gets higher, the noise in the frequency components gets much more amplified, and as the frequency gets lower, the noise in the frequency components gets less and less amplified.
\end{Remark}

\begin{Remark}
From Proposition \ref{Prop Reg method}, we can see that the family of function $(\varphi_s)_{s>0}$ defined by \eqref{def of our molllif kernel} allows to define a family of regularization methods, each regularization method being determined by the choice of the free parameter $s>0$. For instance, the choice $s = 1$ means considering a Cauchy convolution kernel while $s=2$ means taking a Gaussian convolution kernel.
\end{Remark}
Let us end this section by the following Lemma which gives rates of convergence of the mollifier operator $M_\alpha$ corresponding to kernel $\varphi$ defined by \eqref{def of our molllif kernel} on Sobolev spaces $H^p(\Rb^n)$ with $p>0$.
\begin{Lemma}
\label{Lemma rate converg mollifier operator}
Let $p>0$ and $M_\alpha$ be the mollifier operator defined by \eqref{def oper Cbeta} with the mollifier kernel $\varphi$ given by \eqref{def of our molllif kernel}. Then
\begin{equation}
\label{speed conv mollifier operator}
\forall f \in H^p(\Rb^n), \quad || f - M_\alpha f||_{L^2} \leq \tau^{\frac{p \wedge s}{s} }\alpha^{p \wedge s} ||f||_{H^p}.
\end{equation}
\end{Lemma}

\begin{Proof}
Let $f \in  H^p(\Rb^n)$.
If $p<s$, using Parseval identity, we have
\begin{eqnarray*}
\norm{ f - M_\alpha f}_{L^2} & = & \norm{ [1 - \sqrt{2 \pi}^n \widehat{\varphi}(\alpha \xi)] \widehat{f}(\xi) }_{L^2} \\
& = & \norm{ [ 1- \exp(-\tau (\alpha |\xi|)^s)]^{p/s} \widehat{f}(\xi) \times [1- \exp(-\tau (\alpha |\xi|)^s)]^{1 - p/s}  }_{L^2} \\
& \leq & \norm{ [ 1- \exp(-\tau (\alpha |\xi|)^s)]^{p/s} \widehat{f}(\xi)}_{L^2}  \\
& \leq & (\tau \alpha ^s)^{p/s} \norm{(|\xi|^s)^{p/s}\widehat{f}(\xi)}_{L^2} \leq \tau^{p/s} \alpha^{p} \norm{f}_{H^p}.
\end{eqnarray*}
If $p\geq s$, $ 
\norm{ f - M_\alpha f}_{L^2}  =  \norm{ [ 1- e^{-\tau (\alpha |\xi|)^s}] \widehat{f}(\xi)}_{L^2}  \leq  \tau \alpha ^s \norm{|\xi|^s \widehat{f}(\xi)}_{L^2} \leq \tau \alpha^{s} \norm{f}_{H^s} \leq \tau  \alpha^{s} \norm{f}_{H^p}.
$
\end{Proof}
%%---------------------------------------------------------------
%%---------------------------------------------------------------
\section{Error estimates}\label{section error estimate}
Henceforth, $\varphi$ denotes the mollifier kernel defined by \eqref{def of our molllif kernel} and  $g^\delta \in L^2(\Rb^n)$ denotes a noisy data satisfying the noise level condition
\begin{equation}
\label{noise level cond on data}
|| g - g^\delta ||_{L^2} \leq \delta,
\end{equation}
where $g = u(\cdot,T)$ is the exact final distribution. Let us introduce the regularized solution $u_\alpha^\delta$ corresponding to the noisy data $g^\delta$ as the solution of equation
\begin{equation}
\label{equation u_alpha delta}
\begin{cases}
\frac{\partial^\gamma u}{\partial t^\gamma}  = \Delta u  & x \in \Rb^n, \,\,\, t \in (0,T)\\
u(x,T) = (M_\alpha g^\delta)(x) & x \in \Rb^n,
\end{cases}
\end{equation}
Equivalently, we can define $u_\alpha^\delta$ in the frequency domain by
\begin{equation}
\label{express u alpha in freq domain delta}
\widehat{u_\alpha^\delta}(\xi,t) = \sqrt{2 \pi}^n \frac{E_{\gamma,1}(-|\xi|^2 t^\gamma)}{E_{\gamma,1}(-|\xi|^2 T^\gamma)} \widehat{\varphi}(\alpha\xi) \widehat{g^\delta}(\xi) , \quad t \in [0,T].
\end{equation}
It is well known that without assuming a smoothness condition on the exact solution $u(\cdot,0)$ (or on the exact data $g$), it is impossible to exhibit a rate of convergence of regularized solution towards the exact solution \cite{schock1985approximate}. Henceforth, we consider the following classical Sobolev smoothness condition:
\begin{equation}
\label{smoothness cond on u(cdot,0)}
u(\cdot,0) \in H^p(\Rb^n),\quad \norm{u(\cdot,0)}_{H^p} \leq E, \quad \text{with} \quad  p>0, \,\,E>0.
\end{equation}
Before presenting the main results of this section, let us state some lemmas which will be useful in the sequel.
\begin{Lemma}
\label{Lemma bound H p+2 Hp}
Let $p\geq 0$, and $v$ be a solution of equation $\frac{\partial^\gamma v}{\partial t^\gamma}  = \Delta v$ on $\Rb^n$. If $v(\cdot,0) \in H^p(\Rb^n)$, then for every $t \in (0,T]$, $v(\cdot,t) \in H^{p+2}(\Rb^n)$ and 
\begin{equation}
\label{bound H p+2 Hp }
|| v(\cdot,t)||_{H^{p+2}} \leq \frac{C}{1 \wedge t^\gamma } || v(\cdot,0)||_{H^{p}}, \quad \text{with} \quad  C = \frac{C_2}{\Gamma(1-\gamma)}.
\end{equation}
\end{Lemma}
\begin{Proof}
The proof follows readily by applying Parseval identity and estimate \eqref{bounds E gamma Xi t } to equation $\widehat{v}(\xi,t) = \widehat{v}(\xi,0) E_{\gamma,1}(-|\xi|^2 t^\gamma)$ which links $v(\cdot,0)$ and $v(\cdot,t)$ in the frequency domain.
\end{Proof}\\
The next lemma illustrates the fact that the Sobolev smoothness condition \eqref{smoothness cond on u(cdot,0)} is nothing but a Hölder source condition.
\begin{Lemma}
\label{Lemma link smoothness cond and holder source condition}
Let $p > 0$, $u$ be the solution of problem \eqref{main equation} and $A$ being the forward diffusion operator defined in \eqref{def operator A pb}. The smoothness condition $u(\cdot,0) \in H^p(\Rb^n)$ is equivalent to the Hölder source condition $u(\cdot,0) = (A^*A)^{p/4} w$ with $w \in L^2(\Rb^n)$, satisfying
\begin{equation}
\label{estima norm w and u}
 \left(\frac{\Gamma(1-\gamma) (1 \wedge T^\gamma)}{C_2} \right)^{p/2} || u(\cdot,0)||_{H^p} \leq ||w ||_{L^2} \leq \left(\frac{\Gamma(1-\gamma) (1 \vee T^\gamma)}{C_1} \right)^{p/2} || u(\cdot,0)||_{H^p}
\end{equation}
\end{Lemma}

\begin{Proof}
For $u(\cdot,0) \in H^p(\Rb^n)$, formally define $w$ in the frequency domain by
$$
\widehat{w}(\xi) = E_{\gamma,1}(-|\xi|^2 T^\gamma)^{-p/2} \widehat{u}(\xi,0), \quad \Longleftrightarrow \quad \widehat{u}(\xi,0) = \left( E_{\gamma,1}(-|\xi|^2 T^\gamma)^2\right)^{p/4} \widehat{w}(\xi)
$$
From \eqref{def operator A pb}, we can verify that the above definition of $w$ from $u(\cdot,0)$ is merely reformulation of the equation $u(\cdot,0) = (A^*A)^{p/4} w$ in the frequency domain. Next, we can check that $w$ is well defined and belongs to $L^2(\Rb^n)$. Finally, estimate \eqref{estima norm w and u} is deduced from \eqref{bounds E gamma Xi t }.
\end{Proof}

\begin{Remark}
From Lemma \ref{Lemma link smoothness cond and holder source condition}, we can deduce that the order optimal convergence rate under smoothness condition \eqref{smoothness cond on u(cdot,0)} is nothing but $C E^{\frac{2}{p+2}} \delta^{\frac{p}{p+2}}$ with $C \geq 1$ independent of $E$ and $\delta$.
\end{Remark}
The next Lemma which generalizes \cite[Lemma 3]{van2020mollification} will be useful in the sequel for establishing Sobolev norm error estimates.
\begin{Lemma}
\label{Lemma yy}
Let $p\geq 0$ and $v$ be a solution of equation $\frac{\partial^\gamma v}{\partial t^\gamma}  = \Delta v$ on $\Rb^n$. If $v(\cdot,0) \in H^p(\Rb^n)$ then
\begin{equation}
\label{key est H l 1 yy}
\forall l \in [0,p], \quad \forall t \in (0,T], \quad \norm{v(\cdot,t)}_{H^{l+2}} \leq \frac{C(\gamma)}{(1 \wedge t^\gamma)} \norm{v(\cdot,0)}_{H^p}^{\frac{2+l}{p+2}} \norm{v(\cdot,T)}_{L^2}^{\frac{p-l}{p+2}}.
\end{equation}
where 
$
C(\gamma) =  (1 \vee T^\gamma)C_2/C_1$.
Moreover, 
\begin{equation}
\label{key est H l 0 yy}
\forall l \in [0,p], \quad \forall t \in [0,T], \quad  \norm{v(\cdot,t)}_{H^l} \leq \bar{C}(\gamma)  \norm{v(\cdot,0)}_{H^p}^{\frac{2+l}{p+2}} \norm{v(\cdot,T)}_{L^2}^{\frac{p-l}{p+2}},
\end{equation}
where $\bar{C}(\gamma) = \Gamma(1-\gamma)(1 \vee T^\gamma) /C_1$.
\end{Lemma}
\begin{Proof}
Let $l \in [0,p]$ and $v$ be a solution of equation $\frac{\partial^\gamma v}{\partial t^\gamma}  = \Delta v$ with $v(\cdot,0) \in H^p(\Rb^n)$. For $t \in (0,T]$, using Hölder inequality, we have
\begin{eqnarray}
\label{eq 00xx11}
\norm{v(\cdot,t)}_{H^{l+2}}^2 & = & \int_{\Rb^n} (1+ |\xi|^2)^{l+2} \mod{E_{\gamma,1}(-|\xi|^2 t^\gamma) \widehat{v}(\xi,0)}^2 \ud \xi \nonumber \\
& \leq & \frac{C_\gamma^2}{(1 \wedge t^\gamma)^2} \int_{\Rb^n} (1+ |\xi|^2)^{l} |\widehat{v}(\xi,0)|^2 \ud \xi  \quad \text{with} \quad C_\gamma =  \frac{C_2}{\Gamma(1-\gamma)} \quad\text{using} \quad \eqref{bounds E gamma Xi t }\nonumber \\
& = & \frac{C_\gamma^2}{(1 \wedge t^\gamma)^2} \int_{\Rb^n} \left( (1+ |\xi|^2)^{\frac{p(l+2)}{p+2}} |\widehat{v}(\xi,0)|^{\frac{2(l+2)}{p+2}} \right) \left( (1+ |\xi|^2)^{\frac{2(l-p)}{p+2}} |\widehat{v}(\xi,0)|^{\frac{2(p-l)}{p+2}} \right) \ud \xi \nonumber\\
& \leq & \frac{C_\gamma^2}{(1 \wedge t^\gamma)^2}  \left( \int_{\Rb^n} (1+ |\xi|^2)^{p} |\widehat{v}(\xi,0)|^{2}  \ud \xi \right)^{\frac{l+2}{p+2}} \left( \int_{\Rb^n} (1+ |\xi|^2)^{-2} |\widehat{v}(\xi,0)|^{2}  \ud \xi \right)^{\frac{p-l}{p+2}} \\
& \leq &  \left(\frac{C_2}{C_1} \frac{1 \vee T^\gamma}{(1 \wedge t^\gamma)}\right)^2 \norm{v(\cdot,0)}_{H^p}^{\frac{2(l+2)}{p+2}} \left( \int_{\Rb^n} \mod{ E_{\gamma,1}(-|\xi|^2 T^\gamma) \widehat{v}(\xi,0)}^{2}  \ud \xi \right)^{\frac{p-l}{p+2}}  \quad\text{using} \quad \eqref{bounds E gamma Xi t }\nonumber \\
& = &\left(\frac{C_2}{C_1} \frac{1 \vee T^\gamma}{(1 \wedge t^\gamma)}\right)^2 \norm{v(\cdot,0)}_{H^p}^{\frac{2(l+2)}{p+2}} \norm{v(\cdot,T)}_{L^2}^{\frac{2(p-l)} {p+2}},\nonumber
\end{eqnarray}
whence \eqref{key est H l 1 yy}.
On the other hand, for every $t \in [0,T]$ and $l \in [0,p]$, we have
\begin{eqnarray*}
\label{eq 00xx22}
\norm{v(\cdot,t)}_{H^{l}}^2 & = & \int_{\Rb^n} (1+ |\xi|^2)^{l} \mod{E_{\gamma,1}(-|\xi|^2 t^\gamma) \widehat{v}(\xi,0)}^2 \ud \xi \nonumber \\
& \leq & \int_{\Rb^n} (1+ |\xi|^2)^{l}  \mod{\widehat{v}(\xi,0)}^2 \ud \xi \nonumber \\
& \leq &  \left( \int_{\Rb^n} (1+ |\xi|^2)^{p} |\widehat{v}(\xi,0)|^{2}  \ud \xi \right)^{\frac{l+2}{p+2}} \left( \int_{\Rb^n} (1+ |\xi|^2)^{-2} |\widehat{v}(\xi,0)|^{2}  \ud \xi \right)^{\frac{p-l}{p+2}} \nonumber \quad \text{from} \quad \eqref{eq 00xx11}\\
& \leq &  \left(\frac{(1 \vee T^\gamma) \Gamma(1-\gamma)}{C_1}\right)^2 \norm{v(\cdot,0)}_{H^p}^{\frac{2(l+2)}{p+2}} \left( \int_{\Rb^n} \mod{ E_{\gamma,1}(-|\xi|^2 T^\gamma) \widehat{v}(\xi,0)}^{2}  \ud \xi \right)^{\frac{p-l}{p+2}}  \quad\text{using} \quad \eqref{bounds E gamma Xi t }\nonumber \\
& = &\left(\frac{(1 \vee T^\gamma) \Gamma(1-\gamma)}{C_1}\right)^2\norm{v(\cdot,0)}_{H^p}^{\frac{2(l+2)}{p+2}} \norm{v(\cdot,T)}_{L^2}^{\frac{2(p-l)}{p+2}},
\end{eqnarray*}
whence \eqref{key est H l 0 yy}.
\end{Proof}\\
We are ready to state the first main result of this section, which exhibits convergence rates of the reconstruction error $u(\cdot,0) - u_\alpha^\delta(\cdot,0)$ in Sobolev spaces $H^l(\Rb^n)$ with $l \geq 0$.
\begin{Theorem}
\label{Theorem 1}
Assume that the solution $u$ of problem \eqref{main equation} satisfies the smoothness condition \eqref{smoothness cond on u(cdot,0)}. Consider a noisy approximation $g^\delta$ satisfying \eqref{noise level cond on data} and let $u_\alpha^\delta$ be the regularized solution defined by \eqref{express u alpha in freq domain delta}. Then for the a-priori selection rule $\alpha(\delta) = (\delta/E)^{\frac{1}{p+2}}$, we have
\begin{equation}
\label{first conv rate}
\forall l \in [0,p], \quad \text{such that} \quad p-l \leq s, \quad \norm{u(\cdot,0)- u_{\alpha(\delta)}^\delta(\cdot,0)}_{H^l} \leq C E^{\frac{2+l}{p+2}} \delta^{\frac{p-l}{p+2}},
\end{equation}
where $C$ is a constant independent of $\delta$ and $E$.
\end{Theorem}
\begin{Proof}
Using the Parseval identity, we have
\begin{eqnarray}
\label{est reg error}
\norm{u(\cdot,0)- u_\alpha(\cdot,0)}_{H^l} = \norm{\mod{1 - \sqrt{2 \pi}^n\widehat{\varphi}(\alpha \xi)} (1+|\xi|^2)^{l/2} \mod{\widehat{u}(\xi,0)} }_{L^2}  =   \norm{\tilde{u} - M_\alpha \tilde{u}}_{L^2},
\end{eqnarray}
where 
$
\tilde{u} = \Fc^{-1} \left((1+ |\xi|^2)^{l/2} \widehat{u}(\xi,0) \right).
$
Since $u(\cdot,0) \in H^p(\Rb^n)$, $\tilde{u} \in H^{p-l}(\Rb^n)$, then applying \eqref{speed conv mollifier operator} to \eqref{est reg error}, we deduce that if $p - l \leq s$, then
\begin{equation}
\label{est reg erro fin}
\norm{u(\cdot,0)- u_\alpha(\cdot,0)}_{H^l} \leq \tau^{\frac{p-l}{s} }\alpha^{p -l} ||\tilde{u}||_{H^{p-l}} = \tau^{\frac{p-l}{s} }\alpha^{p -l} ||u(\cdot,0)||_{H^{p}} \leq \tau^{\frac{p-l}{s} }\alpha^{p -l} E.
\end{equation}
On the other hand, using \eqref{bounds E gamma Xi t } and \eqref{noise level cond on data}, we have
\begin{eqnarray}
\label{est data erro }
\norm{u_\alpha(\cdot,0)- u_\alpha^\delta(\cdot,0)}_{H^l} & = & \norm{(1+ |\xi|^2)^{l/2}\sqrt{2 \pi}^n \widehat{\varphi}(\alpha \xi) \frac{\widehat{g}(\xi) - \widehat{g^\delta}(\xi)}{E_{\gamma,1}(-|\xi|^2 T^\gamma)}  }_{L^2} \nonumber \\
& \leq & \delta (1 \vee T^\gamma)  \frac{\Gamma(1-\gamma)}{C_1}\norm{(1+ |\xi|^2)^{1 + l/2} \exp(-\tau (\alpha |\xi|)^s ) }_{\infty}.
\end{eqnarray}
But $(1+ |\xi|^2)^{1 + l/2} \exp(-\tau (\alpha |\xi|)^s )  = \left( f_{b,d}(|\xi|^2) \right)^{1+l/2}$ with $d = s/2$ and $b = \frac{\tau \alpha^s}{1+ l/2} \to 0$ as $\alpha \to 0$. Hence applying Lemma \ref{Lemma bound our kernel}, we get that there exists a constant $C$ independent of $\alpha$, such that
\begin{equation}
\label{eq xxrr}
\forall \xi \in \Rb^n, \quad (1+ |\xi|^2)^{1 + l/2} \exp(-\tau (\alpha |\xi|)^s )  \leq (C /\alpha^{2})^{1+l/2}.
\end{equation}
Applying \eqref{eq xxrr} to \eqref{est data erro }, we deduce the existence of a constant $\tilde{C}$ independent of $\alpha$, $\delta$ and $E$ such that
\begin{equation}
\label{est data err fin}
\norm{u_\alpha(\cdot,0)- u_\alpha^\delta(\cdot,0)}_{H^l} \leq \tilde{C} \frac{\delta}{\alpha^{2+l}}.
\end{equation}
Finally from \eqref{est reg erro fin} and \eqref{est data err fin}, we deduce that
\begin{equation}
\label{bound recons error H l}
\norm{u(\cdot,0)- u_\alpha^\delta(\cdot,0)}_{H^l}  \leq  \tau^{\frac{p-l}{s} }\alpha^{p -l} E +  \tilde{C} \frac{\delta}{\alpha^{2+l}}.
\end{equation}
By choosing $\alpha(\delta) = (\delta/E)^{\frac{1}{p+2}}$ in \eqref{bound recons error H l}, we get \eqref{first conv rate}.
\end{Proof}

\begin{Remark}
By considering $l=0$ in Theorem \ref{Theorem 1}, we get that if $p \leq s$, then 
$$
\norm{u(\cdot,0)- u_\alpha^\delta(\cdot,0)}_{L^2} \leq C E^{\frac{2}{p+2}} \delta^{\frac{p}{p+2}}.
$$
Hence, our regularization method is order-optimal under the classical smoothness condition \eqref{smoothness cond on u(cdot,0)}. Notice that the condition $p \leq s$ is not at all restrictive since the parameter $s$ is freely chosen in $\Rb_+^*$.
\end{Remark}

The next theorem shows rate of convergence of the error $u(\cdot,t) - u_\alpha^\delta(\cdot,t)$ when approximating earlier distribution $u(\cdot,t)$, $t \in (0,T)$ using the regularized solution $u_\alpha^\delta$.
\begin{Theorem}
\label{Theorem 2}
Consider the setting of Theorem \ref{Theorem 1}. By considering $\alpha(\delta) = (\delta/E)^{\frac{1}{p+2}}$, we have
\begin{equation}
\label{second conv rate}
\forall t \in (0,T], \,\,\, \forall l \in [0,p+2], \,\,\, \text{s. t.} \,\,\, p+2-l \leq s, \,\,\, \norm{u(\cdot,t)- u_{\alpha(\delta)}^\delta(\cdot,t)}_{H^l} \leq \frac{C}{t^\gamma} E^{\frac{l}{p+2}} \delta^{1-\frac{l}{p+2}},
\end{equation}
where $C$ is a constant independent of $\delta$ and $E$. Moreover, 
\begin{equation}
\label{sec conv rate}
\forall t \in (0,T], \,\,\, \forall l \in [0,p], \,\,\, \text{s. t.} \,\,\, p-l \leq s, \,\,\, \norm{u(\cdot,t)- u_{\alpha(\delta)}^\delta(\cdot,t)}_{H^l} \leq \tilde{C} E^{\frac{2+l}{p+2}} \delta^{\frac{p-l}{p+2}},
\end{equation}
where $\tilde{C}$ is a constant independent of $\delta$, $E$ and $t$.
\end{Theorem}

\begin{Proof}
By noticing that $\mod{u(\xi,t)- u_\alpha^\delta(\xi,t)} = E_{\gamma,1}(-|\xi|^2 t^\gamma) \mod{u(\xi,0)- u_\alpha^\delta(\xi,0)} \leq \mod{u(\xi,0)- u_\alpha^\delta(\xi,0)}$, we can deduce \eqref{sec conv rate} from \eqref{first conv rate}. Let us prove \eqref{second conv rate}.
We have
\begin{eqnarray*}
\label{eqn a}
\norm{u(\cdot,t)- u_\alpha(\cdot,t)}_{H^l} = \norm{\mod{1 - \sqrt{2 \pi}^n\widehat{\varphi}(\alpha \xi)} (1+|\xi|^2)^{l/2} \mod{\widehat{u}(\xi,t)} }_{L^2} = \norm{\tilde{u} - M_\alpha \tilde{u}}_{L^2},
\end{eqnarray*}
where 
$
\tilde{u} = \Fc^{-1} \left((1+ |\xi|^2)^{l/2} \widehat{u}(\xi,t) \right).
$
Since $u(\cdot,0) \in H^p(\Rb^n)$, from Lemma \ref{Lemma bound H p+2 Hp}, we get $u(\cdot,t) \in H^{p+2}(\Rb^n)$ which implies that $\tilde{u} \in H^{p+2-l}(\Rb^n)$, then applying \eqref{speed conv mollifier operator} and \eqref{bound H p+2 Hp }, we get that if $p +2 - l \leq s$, then
\begin{eqnarray}
\label{est reg erro fin bis}
\norm{u(\cdot,t)- u_\alpha(\cdot,t)}_{H^l} \leq \tau^{\frac{p+2-l}{s} }\alpha^{p +2-l} ||\tilde{u}||_{H^{p+2-l}}& = &\tau^{\frac{p+2-l}{s} }\alpha^{p+2 -l} ||u(\cdot,t)||_{H^{p+2}} \nonumber \\
& \leq & \frac{C_2\,\tau^{\frac{p+2-l}{s} }}{(1 \wedge t^\gamma)\Gamma(1-\gamma)}\alpha^{p +2-l} E.
\end{eqnarray}
On the other hand, using \eqref{bounds fract E gamma Xi t }, we have
\begin{eqnarray}
\label{est data erro bis}
\norm{u_\alpha(\cdot,t)- u_\alpha^\delta(\cdot,t)}_{H^l} & = & \norm{(1+ |\xi|^2)^{l/2}\sqrt{2 \pi}^n \widehat{\varphi}(\alpha \xi)\frac{E_{\gamma,1}(-|\xi|^2 t^\gamma)}{E_{\gamma,1}(-|\xi|^2 T^\gamma)} \mod{ \widehat{g}(\xi) - \widehat{g^\delta}(\xi) } }_{L^2} \nonumber \\
& \leq & \delta \,\frac{C_2}{C_1} \left( \frac{T}{t}\right)^\gamma   \norm{(1+ |\xi|^2)^{l/2} \exp(-\tau (\alpha |\xi|)^s ) }_{\infty}.
\end{eqnarray}
When $l=0$, we can easily see from \eqref{est data erro bis} that 
\begin{equation}
\label{ert}
\norm{u_\alpha(\cdot,t)- u_\alpha^\delta(\cdot,t)}_{H^l} \leq \delta \,\frac{C_2}{C_1} \left( \frac{T}{t}\right)^\gamma .
\end{equation}
When $l>0$, using similar reasoning yielding to \eqref{eq xxrr}, we have
\begin{equation}
\label{eq xxrr yy}
\forall \xi \in \Rb^n, \quad 
(1+ |\xi|^2)^{l/2} \exp(-\tau (\alpha |\xi|)^s )  \leq (C /\alpha^{2})^{l/2} 
\end{equation}
From \eqref{est data erro bis}, \eqref{ert} and \eqref{eq xxrr yy}, we deduce that for every $l\geq 0$ the exists a constant $C$ independent of $\alpha$, $\delta$ and $t$ such that
\begin{equation}
\label{eq rr tt}
\norm{u_\alpha(\cdot,t)- u_\alpha^\delta(\cdot,t)}_{H^l}  \leq \frac{C}{t^\gamma} \frac{\delta}{\alpha^l}
\end{equation}
Lastly from \eqref{est reg erro fin bis} and \eqref{eq rr tt}, we get that if $p +2 -l \leq s$, then
\begin{equation}
\label{eq final t}
\norm{u(\cdot,t)- u_\alpha^\delta(\cdot,t)}_{H^l} \leq \frac{\bar{C}}{1\wedge t^\gamma} \alpha^{p+2-l} E  + \frac{C}{t^\gamma} \frac{\delta}{\alpha^l},
\end{equation}
where $C$ and $\bar{C}$ are constants independent of $\alpha$, $\delta$, $E$ and $t$. For $\alpha(\delta) = (\delta/E)^{\frac{1}{p+2}}$, we get \eqref{second conv rate}.
\end{Proof}

\begin{Remark}
In Theorem \ref{Theorem 2}, in the estimate \eqref{sec conv rate}, the rate is lower than the rate in \eqref{second conv rate}, however, notice that the factor $1/t^\gamma $ in \eqref{second conv rate} blows up as $t$ decreases to $0$. Notice that the rate in \eqref{sec conv rate} cannot be improved without multiplying by a factor which blows up as $t$ goes to $0$.
\end{Remark}

\begin{Remark}
By choosing $l=0$ in Theorem \ref{Theorem 2}, we get the rate
$
\norm{u(\cdot,t)- u_\alpha^\delta(\cdot,t)}_{L^2} \leq \frac{C}{t^\gamma}  \delta,
$
which means we can also recover earlier distribution $u(\cdot,t)$ with the best possible rate $\delta$. 
%We also have the following rate for the case $l=0$, 
%$
%\norm{u(\cdot,t)- u_\alpha^\delta(\cdot,t)}_{L^2} \leq C E^{\frac{2}{p+2}}\delta^{\frac{p}{p+2}},
%$
%where the constant $C$ is independent of $t$ and thereby does not blows up as $t$ goes to $0$.
\end{Remark}

Now let us study error estimates when both the data $g$ and the operator $A$ are only approximately known. Indeed, though the operator $A$ is explicitly known as \eqref{def operator A pb}, in practical implementation, this operator is only approximated given that the Mittag Leffler function can only be approximated though with desired accuracy \cite{podlubnv1999fractional}.

In the sequel we assume that $\psi_h$ is a positive function defined on $\Rb_+ \times (0,T]$ and satisfying
\begin{equation}
\label{cond appro E gamma}
\forall t \in (0,T], \quad \norm{\frac{\psi_h(|\xi|,t) - E_{\gamma,1}(-|\xi|^2t^\gamma)}{E_{\gamma,1}(-|\xi|^2t^\gamma)}}_{\infty} \leq h.
\end{equation}
With the function $\psi_h$, we can approximate operator $A$ by the operator $A_h$ defined by 
\begin{equation}
\label{def approx operator A h}
A_h = \Fc^{-1} \psi_h(|\xi|,T) \Fc \quad \text{i.e.} \quad \widehat{A_h f}(\xi) = \psi_h(|\xi|,T)  \widehat{f}(\xi).
\end{equation}
From \eqref{cond appro E gamma}, given that for every $t>0$, and $\xi \in \Rb^n$, $E_{\gamma,1}(-|\xi|^2t^\gamma) \in (0,1]$, we can deduce that $||| A -  A_h ||| \leq h$. Let $u_{\alpha}^{\delta,h}$ be the regularized solution defined in the frequency domain by
\begin{equation}
\label{sol approx data and oper}
\begin{cases}
\vspace{0.2cm}
\widehat{u_{\alpha}^{\delta,h}}(\xi,0) = \sqrt{2 \pi}^{n}\widehat{\varphi}(\alpha \xi)\frac{\widehat{g^\delta}(\xi)}{\psi_h(|\xi|,T)} \\
\widehat{u_{\alpha}^{\delta,h}}(\xi,t) = \sqrt{2 \pi}^{n}\widehat{\varphi}(\alpha \xi)\frac{\psi_h(|\xi|,t) }{\psi_h(|\xi|,T)}\widehat{g^\delta}(\xi) & \text{for} \quad t \in (0,T].
\end{cases}
\end{equation}

The next theorem provides error estimates in the approximation of the initial distribution $u(\cdot,0)$ under the practical setting where both the data $g=u(\cdot,T)$ and the forward diffusion operator are only approximately known.
\begin{Theorem}
\label{Theorem approx data and operator}
Consider the setting of Theorem \ref{Theorem 1}. Assume that $h \leq 1/2$, let $\psi_h$ be a function satisfying \eqref{cond appro E gamma} and $u_{\alpha}^{\delta,h}$ be the approximate solution defined in \eqref{sol approx data and oper}. Then for the a-priori selection rule $\alpha(\delta,h) = (h + \delta/E)^{\frac{1}{p+2}}$, the following convergence rate holds:
\begin{equation}
\label{conv rate approx data and oper}
\forall l \in [0,p], \,\,\, \text{such that} \,\,\, p-l \leq s, \quad \norm{u(\cdot,0) - u_{\alpha(\delta,h)}^{\delta,h}(\cdot,0)}_{H^l} \leq C E^{\frac{2+l}{p+2}} (\delta +h \,E )^{\frac{p-l}{p+2}},
\end{equation}
where $C$ is a constant independent of $\delta$, $h$ and $E$.
\end{Theorem}
\begin{Proof}
For simplicity of notation, we set $\psi(|\xi|,t) = E_{\gamma,1}(-|\xi|^2 t^\gamma)$ for all $t \in [0,T]$. 
For every $t \in (0,T]$ and $\xi \in \Rb^n$ we have
\begin{eqnarray}
\label{eqxxx}
\frac{E_{\gamma,1}(-|\xi|^2t^\gamma)}{\psi_h(|\xi|,t)}  \leq  \frac{\psi(|\xi|,t)}{\psi(|\xi|,t) - |\psi(|\xi|,t) - \psi_h(|\xi|,t)|} 
 =  \frac{1}{1 - \frac{|\psi(|\xi|,t) - \psi_h(|\xi|,t)|}{\psi(|\xi|,t)}}  
 \leq  \frac{1}{1 - h} 
\leq  1 + 2 h \leq 2 .
\end{eqnarray}
The first inequality in \eqref{eqxxx} is due to the fact that $\psi(|\xi|,t) \leq |\psi(|\xi|,t) - \psi_h(|\xi|,t)| + \psi_h(|\xi|,t)$, the second inequality comes from \eqref{cond appro E gamma} and the last two inequalities are due to the fact that $h \leq 1/2$.
%Let us define the function $u_{\alpha}^h$ in the frequency domain as
%$$
%\widehat{u_{\alpha}^{h}}(\xi) = \sqrt{2 \pi}^{n}\widehat{\varphi}(\alpha \xi)\frac{ \widehat{g}(\xi)}{\psi_h(|\xi|,T)}.
%$$
Let $l \in [0,p]$ such that $p-l \leq s$, we recall that from \eqref{est reg erro fin}, we have
\begin{eqnarray}
\label{est reg error app data and ope 0}
\norm{u(\cdot,0)- M_\alpha u(\cdot,0)}_{H^l} \leq \tau^{\frac{p-l}{s} } \alpha^{p -l} E.
\end{eqnarray}
%The last inquality is due to the fact that $\tilde{u} \in H^{p-l}(\Rb^n)$ because of \eqref{smoothness cond on u(cdot,0)} together with the fact that $\norm{\tilde{u}}_{H^{p-1}} = \norm{u(\cdot,0)}_{H^p} \leq E$.
%\eqref{est reg erro fin}
By noticing that
\begin{eqnarray}
\label{ttt}
\mod{\widehat{g}(\xi) - \frac{\psi(|\xi|,T)}{\psi_h(|\xi|,T)}\widehat{g^\delta}(\xi)} & \leq  & \mod{\left[ 1 - \frac{\psi(|\xi|,T)}{\psi_h(|\xi|,T)}\right] \widehat{g}(\xi) } + \mod{\frac{\psi(|\xi|,T)}{\psi_h(|\xi|,T)}[\widehat{g}(\xi) - \widehat{g^\delta}(\xi)]} \nonumber \\
& \leq & \mod{\left[ \frac{\psi(|\xi|,T) - \psi_h(|\xi|,T)}{\psi(|\xi|,T)}\right]  \frac{\psi(|\xi|,T)}{\psi_h(|\xi|,T)}\widehat{g}(\xi) } + 2 \mod{\widehat{g}(\xi) - \widehat{g^\delta}(\xi)} \nonumber \\
& \leq & 2 h \mod{\widehat{g}(\xi) } + 2 \mod{\widehat{g}(\xi) - \widehat{g^\delta}(\xi)}  \quad \text{using} \quad \eqref{cond appro E gamma} \quad \text{and} \quad \eqref{eqxxx},
\end{eqnarray}
we deduce that
\begin{eqnarray}
\label{data err prop approx data and opera}
\norm{M_\alpha u(\cdot,0) - u_\alpha^{\delta,h}(\cdot,0)}_{H^l} & = & \norm{\sqrt{ 2 \pi}^n \widehat{\varphi}(\alpha \xi) (1+ |\xi|^2)^{l/2}\left[ \frac{\widehat{g}(\xi)}{\psi(|\xi|,T)} - \frac{\widehat{g^\delta}(\xi)}{\psi_h(|\xi|,T)} \right] }_{L^2} \nonumber \\
& = & \norm{\sqrt{ 2 \pi}^n \widehat{\varphi}(\alpha \xi) \frac{(1+ |\xi|^2)^{l/2}}{\psi(|\xi|,T)}\left[ \widehat{g}(\xi) - \frac{\psi(|\xi|,T)}{\psi_h(|\xi|,T)}\widehat{g^\delta}(\xi) \right] }_{L^2} \nonumber \\
& \leq &  \left(  2 h \norm{g}_{L^2} + 2 \delta \right) \norm{\frac{\exp(-\tau(\alpha \xi)^s) (1+ |\xi|^2)^{l/2}}{E_{\gamma,1}(-|\xi|^2T^\gamma)} }_{\infty} \,\, \text{using} \quad \eqref{ttt} \,\,\, \text{and} \,\,\, \eqref{noise level cond on data} \nonumber \\
& \leq & \left(  2 h E + 2 \delta \right) C_\gamma \norm{(1+ |\xi|^2)^{1 + l/2} \exp(-\tau (\alpha |\xi|)^s ) }_{\infty} \,\, \text{using} \quad \eqref{bounds E gamma Xi t } \,\,\, \text{and} \,\,\, \eqref{smoothness cond on u(cdot,0)} \nonumber \\
& \leq & C \frac{h E +  \delta }{\alpha^{2+l}} \quad \text{using} \quad \eqref{eq xxrr},
\end{eqnarray}
where $C$ is a constant independent of $\alpha$, $\delta$ and $E$.
Finally, from \eqref{est reg error app data and ope 0} and \eqref{data err prop approx data and opera}, we get
$$
\norm{u(\cdot,0) - u_\alpha^{\delta,h}}_{H^l} \leq \tau^{\frac{p-l}{s} } \alpha^{p -l} E  + C E\, \frac{h  +  \delta/E }{\alpha^{2+l}}
$$
from which \eqref{conv rate approx data and oper} follows by choosing $\alpha(\delta,h) = (h + \delta/E)^{\frac{1}{p+2}}$.
\end{Proof}

\begin{Remark}
From Theorem \ref{Theorem approx data and operator}, by choosing $l=0$ in \eqref{conv rate approx data and oper}, we get the rate 
$$
\norm{u(\cdot,0) - u_{\alpha}^{\delta,h}}_{L^2} \leq C E^{\frac{2}{p+2}} (\delta + h \,E)^{\frac{p}{p+2}}.
$$
Hence in the practical setting where both the data and the operator are approximated, we are able to derive order-optimal convergence rates.
\end{Remark}
The next theorem exhibits error estimates when approximating $u(\cdot,t)$ with $t \in (0,T]$.
\begin{Theorem}
\label{Theo conv rates u t data and oper approx}
Consider the setting of Theorem \ref{Theorem approx data and operator}. Then for the a-priori selection rule $\alpha(\delta,h) = (h + \delta/E)^{\frac{1}{p+2}}$, the following convergence rate holds:
\begin{equation}
\label{conv rate approx data and oper t}
\forall t \in (0,T], \,\, \forall l \in [0,p+2], \,\, \text{s. t.} \,\, p+2-l \leq s, \,\, \norm{u(\cdot,t) - u_{\alpha(\delta,h)}^{\delta,h}(\cdot,t)}_{H^l} \leq \frac{C}{t^\gamma} E^{\frac{l}{p+2}} (\delta + h \, E)^{1 - \frac{l}{p+2}},
\end{equation}
where $C$ is a constant independent of $\delta$, $h$, $E$ and $t$.
\end{Theorem}

\begin{Proof} Let $u$ satisfies \eqref{smoothness cond on u(cdot,0)} and $t\in (0,T]$. For $ l \in [0,p+2]$ such that $p +2 - l \leq s$, recall that from \eqref{est reg erro fin bis}, we have
\begin{equation}
\label{Eq xxyy}
\norm{u(\cdot,t)- M_\alpha u(\cdot,t)}_{H^l} \leq  \frac{C_2\,\tau^{\frac{p+2-l}{s} }}{(1 \wedge t^\gamma)\Gamma(1-\gamma)}\alpha^{p +2-l} E.
\end{equation}
On the other hand, using \eqref{cond appro E gamma}, we have
\begin{eqnarray}
\label{ttt 1}
\mod{\widehat{g}(\xi) - \frac{\psi_h(|\xi|,t)}{\psi(|\xi|,t)} \frac{\psi(|\xi|,T)}{\psi_h(|\xi|,T)}\widehat{g^\delta}(\xi)} & \leq  & \mod{\left[ 1 - \frac{\psi_h(|\xi|,t)}{\psi(|\xi|,t)}\right] \widehat{g}(\xi) } + \mod{\frac{\psi_h(|\xi|,t)}{\psi(|\xi|,t)}\left[\widehat{g}(\xi) - \frac{\psi(|\xi|,T)}{\psi_h(|\xi|,T)}\widehat{g^\delta}(\xi)\right]} \nonumber \\
& \leq & \mod{\left[ \frac{\psi(|\xi|,t) - \psi_h(|\xi|,t)}{\psi(|\xi|,t)}\right] \widehat{g^\delta}(\xi)}   + (1+h) \mod{\widehat{g}(\xi) - \frac{\psi(|\xi|,T)}{\psi_h(|\xi|,T)}\widehat{g^\delta}(\xi)} \nonumber \\
& \leq & h \mod{\widehat{g}(\xi) } + (3/2) \left( 2 h \mod{\widehat{g}(\xi) } + 2 \mod{\widehat{g}(\xi) - \widehat{g^\delta}(\xi)} \right) \quad \text{using} \quad \eqref{ttt} \nonumber \\
& = & 4 h \mod{\widehat{g}(\xi) } + 3 \mod{\widehat{g}(\xi) - \widehat{g^\delta}(\xi)} .
\end{eqnarray}
From \eqref{ttt 1}, we deduce that
\begin{eqnarray}
\label{data err prop approx data and opera t}
\norm{M_\alpha u(\cdot,t) - u_\alpha^{\delta,h}(\cdot,t)}_{H^l} & = & \norm{\sqrt{ 2 \pi}^n \widehat{\varphi}(\alpha \xi) (1+ |\xi|^2)^{l/2}\left[ \frac{\psi(|\xi|,t)}{\psi(|\xi|,T)}\widehat{g}(\xi) - \frac{\psi_h(|\xi|,t)}{\psi_h(|\xi|,T)}\widehat{g^\delta}(\xi) \right] }_{L^2} \nonumber \\
& = & \norm{\sqrt{ 2 \pi}^n \widehat{\varphi}(\alpha \xi) (1+ |\xi|^2)^{l/2} \frac{\psi(|\xi|,t)}{\psi(|\xi|,T)} \left[ \widehat{g}(\xi) - \frac{\psi_h(|\xi|,t)}{\psi(|\xi|,t)} \frac{\psi(|\xi|,T)}{\psi_h(|\xi|,T)}\widehat{g^\delta}(\xi) \right] }_{L^2} \nonumber \\
& \leq &  \left(  4 h \norm{g}_{L^2} + 3 \delta \right)\frac{C_2}{C_1} \left( \frac{T}{t}\right)^\gamma \norm{\exp(-\tau(\alpha \xi)^s) (1+ |\xi|^2)^{l/2} }_{\infty} \,\, \text{using} \,\,\, \eqref{bounds fract E gamma Xi t }, \nonumber \\
& \leq & \frac{C}{t^\gamma} \frac{h E +  \delta }{\alpha^{l}} \quad \text{using} \quad \eqref{eq xxrr yy},
\end{eqnarray}
where $C$ is a constant independent of $\alpha$, $\delta$, $E$ and $t$.
Finally, from \eqref{Eq xxyy} and \eqref{data err prop approx data and opera t}, we get
$$
\norm{u(\cdot,t) - u_\alpha^{\delta,h}(\cdot,t)}_{H^l} \leq \frac{C_2\,\tau^{\frac{p+2-l}{s} }}{(1 \wedge t^\gamma)\Gamma(1-\gamma)}\alpha^{p +2-l} E  + \frac{C}{t^\gamma} E\, \frac{h  +  \delta/E }{\alpha^{l}}
$$
from which \eqref{conv rate approx data and oper t} follows by choosing $\alpha(\delta) = (h + \delta/E)^{\frac{1}{p+2}}$.
\end{Proof}

\begin{Remark}
By choosing $l=0$ in Theorem \ref{Theo conv rates u t data and oper approx}, we recover the best possible rate
$$
\norm{u(\cdot,t) - u_{\alpha}^{\delta,h}(\cdot,t)}_{L^2} \leq \frac{C}{t^\gamma}  (\delta + h \,E).
$$
\end{Remark}

Now let us focus on the choice of the regularization parameter $\alpha$ when we don't have precise a-priori information about the smoothness of the sought solution $u(\cdot,0)$.
%%---------------------------------------------------------------
%%---------------------------------------------------------------
\section{A-posteriori parameter choice rule}\label{section par choice rule}
The choice of the regularization parameter is a crucial step for any regularization method. As a matter of fact, no matter the regularization method considered, a bad choice of the regularization parameter results in poor approximate solution. 

Let us consider the following a-posteriori parameter choice rule based on the Morozov principle \cite{mathe2008general}.
\begin{equation}
\label{def beta a posteriori rule}
\alpha(\delta,g^\delta) = \sup \left\lbrace \alpha > 0, \quad \text{s.t.} \quad ||  g^\delta - M_\alpha g^\delta ||_{L^2} < \theta \delta \right\rbrace,
\end{equation} 
where $\theta > 1$ is a free real parameter.
\begin{Proposition}
\label{Prop exist post par choice rule}
Assume that the noise level $\delta$ and the noisy data $g^\delta$ satisfies
\begin{equation}
\label{noise ration cond}
0 < \theta \delta < ||g^\delta||,
\end{equation}
then the parameter $\alpha(\delta,g^\delta)$ expressed in \eqref{def beta a posteriori rule} is well defined and satisfies
\begin{equation}
\label{char beta a posteriori rule}
|| u_{\alpha(\delta,g^\delta)}^\delta(\cdot,T) - g^\delta|| = \theta  \delta \quad \text{i.e.} \quad \norm{ \mod{1- e^{-\tau (\alpha |\xi|)^s} }\widehat{g^\delta}(\xi) }_{L^2} = \theta \delta.
\end{equation}
\end{Proposition}

\begin{Proof}
Let the function $v: \Rb_+^* \to \Rb_+ $ be defined by $v(\alpha) := \norm{ g^\delta - M_\alpha g^\delta}_{L^2}^2$. By Parseval identity, we get that $v(\alpha) = \norm{ \mod{1- e^{-\tau (\alpha |\xi|)^s} }\widehat{g^\delta}(\xi) }_{L^2}^2$. Using the dominated convergence theorem, we can readily check that $\lim_{\alpha \to 0} v(\alpha) = 0$ and $\lim_{\alpha \to \infty} v(\alpha) = \norm{\widehat{g^\delta}(\xi) }_{L^2}^2 $. Moreover, we can also check using derivation under integral sign that $v$ is differentiable with $v'(\alpha) >0$ for all $\alpha>0$ which implies that $v$ is strictly increasing. Hence if \eqref{noise ration cond} is satisfied, then $\theta \delta$ is in the range of the one-to-one function $v$, whence the existence and uniqueness of $\alpha(\delta,g^\delta)$ defined in \eqref{def beta a posteriori rule} which is characterized by \eqref{char beta a posteriori rule}.
\end{Proof}
\begin{Remark}
The condition \eqref{noise ration cond} is quite reasonable as we do not expect to recover reasonable approximate solution if the data is dominated by noise.
\end{Remark}

\begin{Theorem}
\label{Theo rate  posteriori rule 1}
Consider the setting of Theorem \ref{Theorem 1}. Assume that \eqref{noise ration cond} is satisfied and let $\alpha(\delta,g^\delta)$ satisfying \eqref{char beta a posteriori rule}. If $p +2  \leq s$, then the following holds:
\begin{equation}
\label{rate posteriori rule 1}
\forall l \in [0,p], \quad \forall t \in [0,T],, \quad \norm{ u(\cdot,t) - u_{\alpha(\delta,g^\delta)}^\delta(\cdot,t)}_{H^l} \leq C E^{\frac{2+l}{p+2}} \delta^{\frac{p-l}{p+2}}
\end{equation}
where $C$ is a constant independent of $\delta$, $E$ and $t$.
\end{Theorem}

\begin{Proof}
For simplicity of notation, let $\alpha:=\alpha(\delta,g^\delta)$ defined by \eqref{char beta a posteriori rule} and let $g_\alpha = M_\alpha g$ where $g = u(\cdot,T)$ and $g_\alpha^\delta = M_\alpha g^\delta$ . We have
\begin{equation}
\label{eq 001}
\norm{g - g_\alpha} \leq \norm{g - g^\delta} + \norm{g^\delta - g_\alpha^\delta } + \norm{g_\alpha^\delta - g_\alpha} \leq \delta + \theta \delta + \norm{e^{-\tau (\alpha |\xi|)^s} (\widehat{g} - \widehat{g^\delta}(\xi)} \leq (\theta +2)\delta.
\end{equation}
Let $p \geq 0$ such that $p+2 \leq s$, using \eqref{speed conv mollifier operator} and \eqref{bound H p+2 Hp }, we have
\begin{eqnarray}
\label{eq 002}
\theta \delta = \norm{\mod{1-e^{-\tau (\alpha |\xi|)^s} } \widehat{g^\delta}(\xi)} & \leq & \norm{\mod{1-e^{-\tau (\alpha |\xi|)^s} } \widehat{g}(\xi)} + \norm{\mod{1-e^{-\tau (\alpha |\xi|)^s} }\mod{ \widehat{g}(\xi)- \widehat{g^\delta}(\xi)}}  \nonumber \\
& \leq & \tau^{(p+2)/s} \alpha^{p+2} ||g||_{H^{p+2}} + \delta \nonumber  \\
&\leq  & C(p,\gamma) \alpha^{p+2} E + \delta
\end{eqnarray}
From \eqref{eq 002}, we deduce that
\begin{equation}
\label{eq 00002}
(\theta -1) \delta \leq C(p,\gamma) \alpha^{p+2} E \quad  \Rightarrow \quad \frac{1}{\alpha} \leq \left(\frac{C(p,\gamma) }{\theta -1} \right)^{\frac{1}{p+2}} E^{\frac{1}{p+2}} \delta^{-\frac{1}{p+2}}
\end{equation}
By noticing that $\norm{u(\cdot,0) - M_\alpha u(\cdot,0) }_{H^p} \leq \norm{u(\cdot,0) }_{H^p} \leq E$, and applying \eqref{key est H l 0 yy} of Lemma \ref{Lemma yy} together with \eqref{eq 001}, we get
\begin{equation}
\label{eq 003}
\forall l \in [0,p], \quad \forall t \in [0,T], \quad \norm{u(\cdot,t) - u_\alpha(\cdot,t)}_{H^l} \leq C(\gamma,p,\theta) E^{\frac{2+l}{p+2}}\delta^{\frac{p-l}{p+2}},
\end{equation}
where $u_\alpha(\cdot,t) = M_\alpha u(\cdot,t)$ is the solution of equation \eqref{equation u_alpha} and $C(\gamma,p,\theta) = \bar{C}(\gamma) (\theta-1)^{\frac{p-l}{p+2}}$.
On the other hand, from \eqref{est data err fin} and \eqref{eq 00002}, we have that for every $l \in [0,p]$ and $t \in [0,T]$,
\begin{eqnarray}
\label{eq 004}
 \norm{u_\alpha(\cdot,t)- u_\alpha^\delta(\cdot,t)}_{H^l}  \leq \norm{u_\alpha(\cdot,0)- u_\alpha^\delta(\cdot,0)}_{H^l}  
\leq \tilde{C} \frac{\delta}{\alpha^{2+l}} 
\leq C E^{\frac{2+l}{p+2}} \delta^{ \frac{p-l}{p+2}},
\end{eqnarray}
where $C$ is a constant independent of $\delta$, $E$ and $t$. Estimate \eqref{rate posteriori rule 1} follows from \eqref{eq 003} and \eqref{eq 004}.
\end{Proof}

\begin{Remark}
In Theorem \ref{Theo rate  posteriori rule 1}, by choosing $l=0$ and $t = 0$ in \eqref{rate posteriori rule 1}, we obtain the rate 
$$
\norm{ u(\cdot,t) - u_{\alpha(\delta,g^\delta)}^\delta(\cdot,t)}_{L^2} \leq C E^{\frac{2}{p+2}} \delta^{\frac{p}{p+2}},
$$
which means that the parameter choice rule \eqref{def beta a posteriori rule} leads to order-optimal rate.
\end{Remark}

Lastly, the next theorem exhibits rates of convergence when approximating early distribution $u(\cdot,t)$ with $t \in (0,T)$ with the parameter choice rule \eqref{def beta a posteriori rule}.
\begin{Theorem}
\label{Theorem  rate posterori rule 2}
Consider the setting of Theorem \ref{Theo rate  posteriori rule 1}. The following estimate holds
\begin{equation}
\label{rate posteriori 2}
\forall t \in (0,T], \,\, \forall l \in [0,p+2], \,\, \text{s. t.} \,\,\, p+2-l \leq s, \,\, \norm{u(\cdot,t) - u_{\alpha(\delta,g^\delta)}^\delta(\cdot,t)}_{H^l} \leq \frac{C}{t^\gamma} E^{\frac{l}{p+2}} \delta^{1 - \frac{l}{p+2}},
\end{equation}
where $C$ is a constant independent of $\delta$, $E$ and $t$.
\end{Theorem}

\begin{Proof}
Let $t \in (0,T]$, $l \in [0, p] $ and $\alpha := \alpha(\delta,g^\delta)$ defined by \eqref{char beta a posteriori rule}. By applying \eqref{key est H l 1 yy}, we get 
\begin{eqnarray}
\label{eq ttyyww}
\norm{u(\cdot,t) - u_\alpha(\cdot,t) }_{H^{l+2}} & \leq & \frac{C(\gamma)}{1\wedge t^\gamma} \norm{u(\cdot,0) - u_\alpha(\cdot,0) }_{H^p}^{\frac{l+2}{p+2}} \norm{u(\cdot,T) - u_\alpha(\cdot,T) }_{L^2}^{\frac{p-l}{p+2}} \nonumber \\
& \leq & \frac{C(\gamma)}{1\wedge t^\gamma} E^{\frac{l+2}{p+2}} \norm{g - g_\alpha }_{L^2}^{\frac{p-l}{p+2}} \nonumber \\
& \leq & \frac{C(\gamma)}{1\wedge t^\gamma} (\theta +2)^{\frac{p-l}{p+2}} E^{\frac{l+2}{p+2}} \delta^{\frac{p-l}{p+2}} \quad \text{from} \quad \eqref{eq 001}.
\end{eqnarray}
From \eqref{eq ttyyww}, we deduce that for all $l \in [0,p+2]$,
\begin{equation}
\label{eq 008}
\norm{u(\cdot,t) - u_\alpha(\cdot,t) }_{H^{l}} \leq \frac{C}{1\wedge t^\gamma} E^{\frac{l}{p+2}} \delta^{\frac{p+2-l}{p+2}}.
\end{equation}
On the other hand, from \eqref{eq rr tt}, we have
\begin{eqnarray}
\label{eq 009}
\norm{u_\alpha(\cdot,t)- u_\alpha^\delta(\cdot,t)}_{H^l} 
 \leq  \frac{C}{t^\gamma} \frac{\delta}{\alpha^l} \leq  \frac{C}{t^\gamma} \left(\frac{C(p,\gamma) }{\theta -1} \right)^{\frac{l}{p+2}} E^{\frac{l}{p+2}} \delta^{1 - \frac{l}{p+2}} \quad \text{using} \quad \eqref{eq 00002}.
\end{eqnarray}
Estimate \eqref{rate posteriori 2} follows readily from \eqref{eq 008} and \eqref{eq 009}.
\end{Proof}

\begin{Remark}
In Theorem \ref{Theorem  rate posterori rule 2}, taking $l=0$, we get the best possible rate 
$
\norm{u(\cdot,t) - u_{\alpha(\delta,g^\delta)}^\delta(\cdot,t)}_{L^2} \leq C \frac{\delta}{t^\gamma} 
$
for the posteriori parameter choice rule \eqref{def beta a posteriori rule}.
\end{Remark}
Let us end this section with the following algorithm for approximating the regularization parameter $\alpha(\delta,g^\delta)$.

\begin{algorithm}%[h !]
\begin{center}
\begin{algorithmic}[1]
\State Set $\alpha_0 \gg 1$ and $q \in (0,1)$
\State Set $\alpha = \alpha_0$ (initial guess)
\While{$|| \mod{1- e^{-\tau (\alpha |\xi|)^s} }\widehat{g^\delta}(\xi) || \,>\, \theta \delta $}
\State $\alpha = q \times \alpha$
\EndWhile
\end{algorithmic}
\end{center}
\caption{}
\label{Algo alpha }
\end{algorithm}
%%-----------------------------------------------------------
%%-----------------------------------------------------------
\section{Numerical experiments}\label{section numerical experiments}
In order to illustrate the effectiveness of our regularization approach, we consider four numerical examples in two-dimension space where we invariably set $T = 1$ and $\gamma = 0.8$.

\textbf{Example 1}: $u(x,0) = e^{-x_1^2 - x_2^2}$.

\textbf{Example 2}: $u(x,0) = e^{-\mod{x_1} - \mod{x_2}}$.

\textbf{Example 3}: $u(x,0) = v(x_1)v(x_2)$ where $v$ is triangle impulse set as $v(\lambda)=
\begin{cases}
1 + \lambda/3 & \text{if} \,\, \lambda \in [-3,0]\\
1 - \lambda/3 & \text{if} \,\, \lambda \in (0,3] \\
0 & \text{otherwise}.
\end{cases}
$

\textbf{Example 4}:  $
u(x,0) =
\begin{cases}
1  & \text{if} \,\, x \in [-5,5]^2\\
0 & \text{otherwise}.
\end{cases}
$

Notice that, in Examples 1 2 and 3, $u$ is continuous on the contrary to Example 4. Moreover, in Example 1, $u(\cdot,0) \in H^{p}(\Rb^2)$ for all $p>0$; in Example 2, $u(\cdot,0) \in H^{p}(\Rb^2)$ for $p<1$; in Example 3, $u(\cdot,0) \in H^{1}(\Rb^2)$; and in Example 4, $u(\cdot,0) \in H^{p}(\Rb^2)$ for $p <1/2$.

In the four examples, the support of $u(\cdot,0)$ is considered $[-L,L]^2$ which is uniformly discretized as $(x_1(i),x_2(j))$ where $x_1(i)=x_2(i) = -L + (i-0.5)\kappa$ with $\kappa = 2L/N$ and $i,j=1,...,N$. In all the simulations, we set $L=10$ and $N =256$.

The noisy data $g^\delta$ is generated as 
$
g^\delta(x) = u(x,T) + \eta\, \epsilon(x),
$
 where $\epsilon(x)$ is a random number drawn from the standard normal distribution and $\eta$ is a parameter allowing to control to amount to noise added to the data. The noise level $\delta$ is nothing but $ \eta  E ||\epsilon||_2 $ and the percentage of noise $perc\_noise$ is nothing but 
 $$ 
 perc\_noise =  \frac{100\times \eta  E ||\epsilon||_2}{||u(\cdot,T)||_2} \,\,\%.
 $$
For the mollifier kernel $\varphi$ defined in \eqref{def of our molllif kernel}, we invariably choose $s=4$ and $\tau = 1/2$.
The Fourier transform and inverse Fourier transform involved in the computation of the reconstructed solution $u_\alpha^\delta$ are quite rapidly evaluated with the numerical procedure from \cite{bailey1994fast} based on fast Fourier transform (FFT) algorithm. From the Shannon-Nyquist principle, we set the frequency domain to $[-\Omega,\Omega]^2$ with $\Omega = \pi N/2L$.

\begin{figure}[h!]
\begin{center}
\includegraphics[scale=0.4]{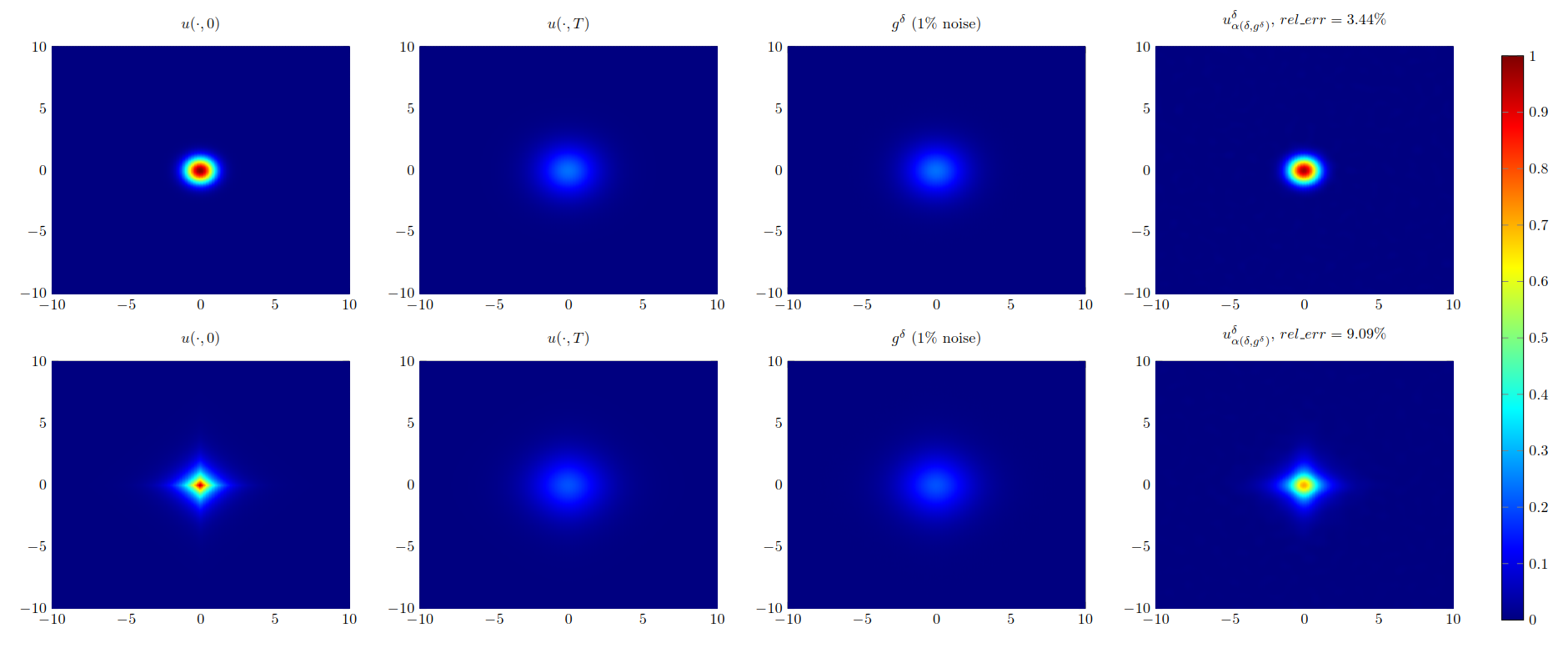} 
\end{center}
\caption{Illustration reconstructed solution $u_{\alpha(\delta,g^\delta)}^\delta$ for $perc\_noise = 1\%$ in Example 1 (first row) and Example 2 (second row).}
\label{Fig comp sol 1 and 2 percent noise level}
\end{figure}

\begin{figure}[h!]
\begin{center}
\includegraphics[scale=0.4]{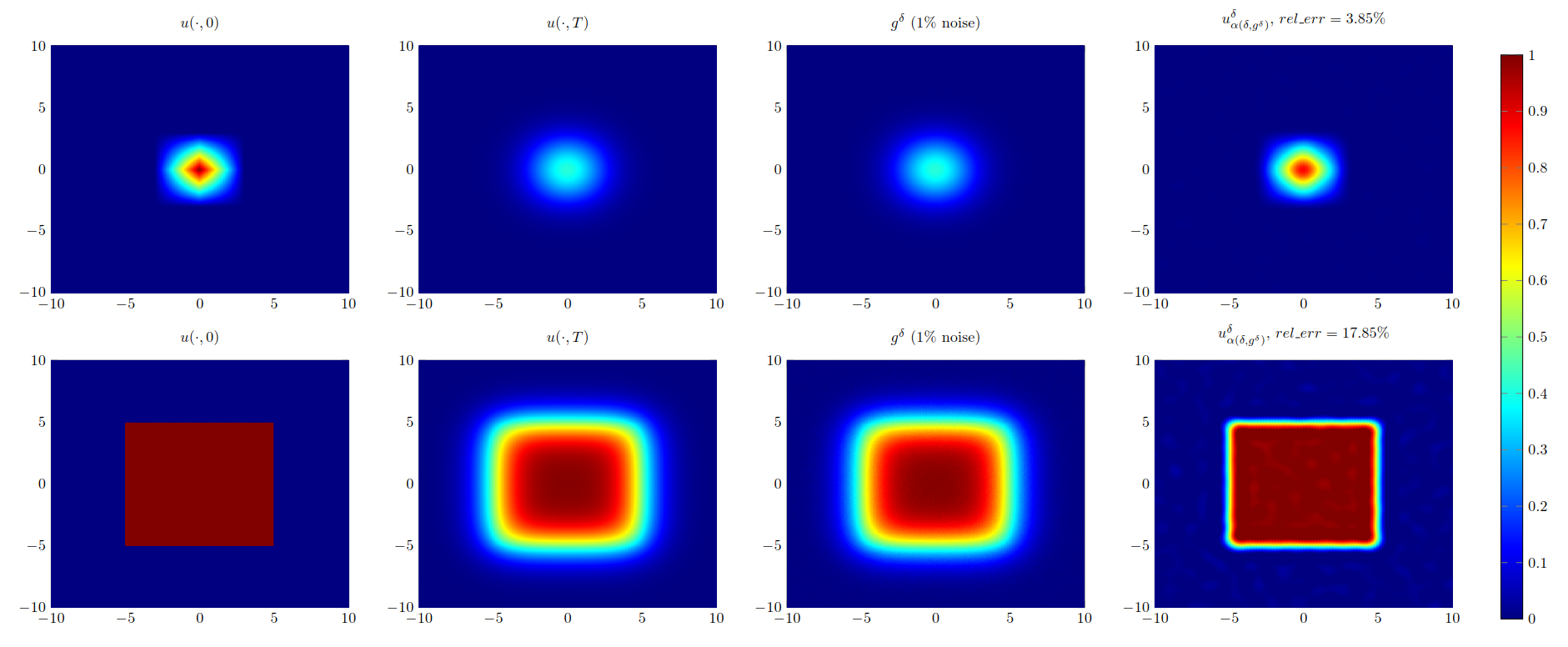} 
\end{center}
\caption{Illustration reconstructed solution $u_{\alpha(\delta,g^\delta)}^\delta$ for $perc\_noise = 1\%$ in Example 3 (first row) and Example 4 (second row).}
\label{Fig comp sol 3 and 4 percent noise level}
\end{figure}

For all the results in this section, we consider the parameter choice rule $\alpha(\delta,g^\delta)$ defined by \eqref{def beta a posteriori rule} with $\theta =1.01$. We compute this parameter via Algorithm \ref{Algo alpha } with $q=0.99$ and $\alpha_0 = 10$. On Figures \ref{Fig comp sol 1 and 2 percent noise level} and \ref{Fig comp sol 3 and 4 percent noise level}, we illustrate the approximate solution for $perc\_noise = 1 \%$ in each example.

 The reconstructed solution $u_{\alpha(\delta,g^\delta)}^\delta$ is computed via formula \eqref{express u alpha in freq domain delta} (with $t=0$) followed by inverse Fourier transform. For each approximate solution $u_{\alpha(\delta,g^\delta)}^\delta$, we assess the relative error
$$
rel\_err = \frac{|| u_{\alpha(\delta,g^\delta)}^\delta- u(\cdot,0)||_2}{||u(\cdot,0)||_2}.
$$

\begin{figure}[h!]
\begin{center}
\includegraphics[scale=0.45]{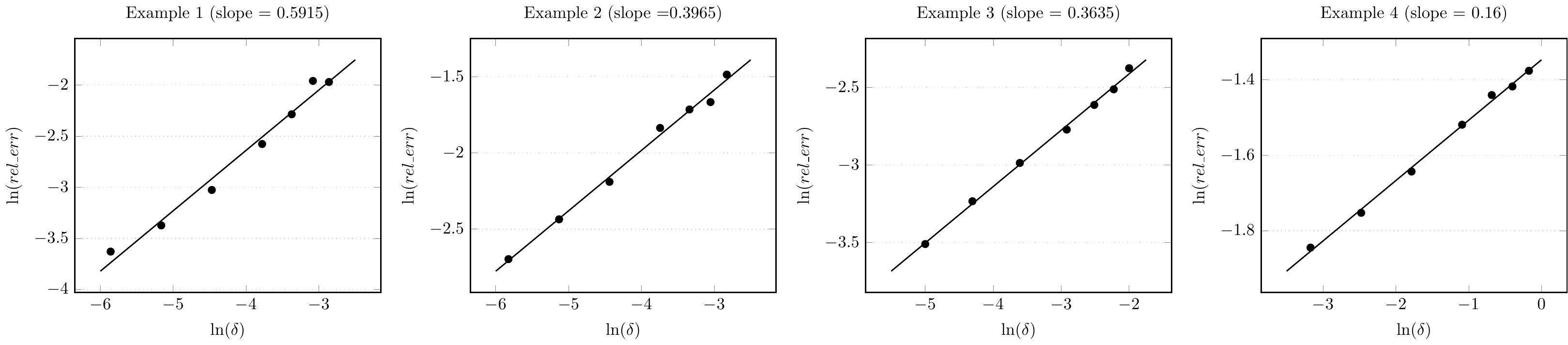} 
\end{center}
\caption{Illustration of numerical rates of convergence for the a-posteriori rule \eqref{def beta a posteriori rule}.}
\label{Fig ill num conv rate}
\end{figure}

In order to numerically confirm the theoretical rates of the reconstruction error, on Figure \ref{Fig ill num conv rate}, we plot $ln(rel\_err)$ versus $ln(\delta)$ for various values of $\delta$. We recall that if the reconstruction error is of order $\mathcal{O}\left(\delta^r \right)$ as $\delta \to 0$, then the curve $(\ln(\delta),\ln(rel\_err))$ should exhibit a line shape with slope equal to $r$. From Figure \ref{Fig ill num conv rate}, we can see that for each example, the plot clearly exhibits a line shape, confirming thus the power rate in the reconstruction error. Though the numerical order $r_{num}$ observed is different from the theoretical rate $r = p/(p+2)$, however, the numerical order does increase as the solution gets smoother. In fact, the higher order ($0.5915$) is achieved for example 1 which corresponds to the smoothest case; next in example 2 and 3, where the solution $u(\cdot,0)$ has approximated the same regularity, the numerical orders in this cases are closed; lastly, the lowest numerical order ($0.16$) is achieved in Example 4 where the solution is the least regular.

\begin{table}[h!]
\begin{center}
\includegraphics[scale=0.9]{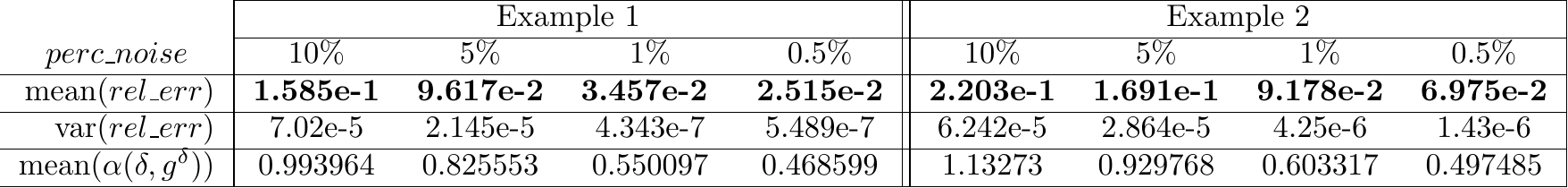}
\end{center}
\caption{Summary Monte Carlo simulation for Example 1 and 2 with 200 sample size.}
\label{Table MC exple 1,2}
\end{table}

\begin{table}[h!]
\begin{center}
\includegraphics[scale=0.9]{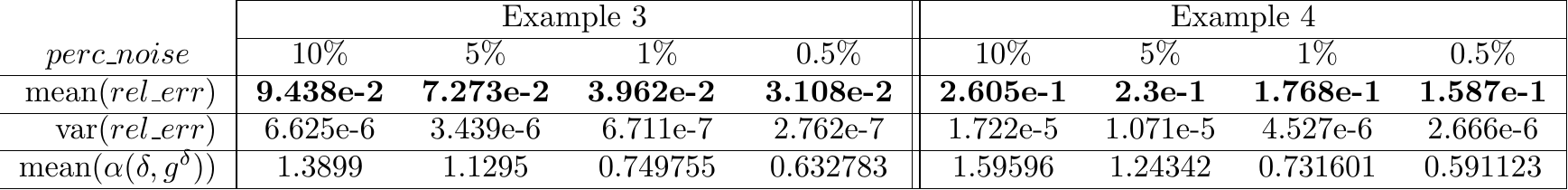} 
\end{center}
\caption{Summary Monte Carlo simulation for Example 3 and 4 with 200 sample size.}
\label{Table MC exple 3,4}
\end{table}

Finally, in order to assess the numerical stability and convergence of our method, we run a Monte Carlo simulations of $200$ replications of noise term for each example with various percentage of noise. The results are summarized in Tables \ref{Table MC exple 1,2} and \ref{Table MC exple 3,4}.  From Tables \ref{Table MC exple 1,2} and \ref{Table MC exple 3,4}, we can see that :
\begin{itemize}
\item The reconstruction error and the regularization parameter $\alpha(\delta,g^\delta)$ both decrease as the noise level decreases. This indicates the numerical convergence of the regularization method coupled with the parameter choice rule \eqref{def beta a posteriori rule}.

\item The variance of the reconstruction error in all cases is not larger than $10^{-4}$. This is a significant indicator of the high numerical stability of the regularization technique coupled with parameter choice rule \eqref{def beta a posteriori rule}.

\item Lastly, the reconstruction error are much smaller in Example 1 and 3 compared to Examples 2 and 4. This observation confirms the following known fact: the smoother the sought solution, the better the reconstruction is expected. Indeed, the solution in Examples 1 and 3 are much smoother than the solution in Examples 2 and 4.
\end{itemize}
%%---------------------------------------------------------------------
\section{Conclusion}	
In this paper, we focused on the regularization of final value time-fractional diffusion equation on unbounded domain. We presented a broad class of regularization methods which encompasses some regularization methods appearing distinctively in literature for this problem. We proposed a simple regularization method which smoothly regularizes the problem without truncating high frequency components. We proved order-optimality of our regularization approach in the practical setting where not only the data are approximated but also the forward diffusion operator. We also showed that with our regularization approach, we can also approximate early distribution $u(\cdot,t)$ with $t \in (0,T)$ with the best possible rate. For full applicability of the method, based on the Morozov principle, we provided a parameter choice rule leading to order-optimal rates. Finally, the soundness and efficiency of the regularization technique couples with the parameter choice rule prescribed is illustrated through some numerical simulations in two-dimensional space.
%\clearpage
\section*{Appendix}

\textbf{Proof of Lemma \ref{Lemma bound our kernel}:}
Let the function $f(x) = (1+ x) e^{- b x^d}$ with $b,d \in \Rb_+^*$.
The function $f$ is continuous on $[0,+\infty )$ with $f(0) = 1$ and $ \lim_{x\to +\infty} f(x) = 0$. Moreover, $f$ is smooth on $(0,\infty)$, hence, $f$ admits a global maximum on $\Rb_+$ attained at a critical point.
Now, notice that $f'(x) = p(x)e^{- b x^d}$ with $p(x) = 1 - b d x^{d-1}(1+x)$. Hence $f'(x) = 0$ implies that $p(x) = 0$. 
But 
\begin{equation}
\label{ee 000}
p(x) = 0 \Rightarrow x^{d-1}(1+x) = \frac{1}{b d} \to + \infty \quad \text{as} \quad b \downarrow 0.
\end{equation}

If $d \geq 1$, we can verify that $p$ has a unique root $\bar{x}$ on $ \Rb_+$ which maximizes $f$. Moreover, $\bar{x} \to +\infty $ as $b \downarrow 0$. Since $\bar{x} \to +\infty $ as $b \downarrow 0$, from \eqref{ee 000}, we deduce that $ \bar{x}^d  \sim 1/bd$. Since $p(\bar{x}) = 0$, then, $(1 +  \bar{x}) = (1/bd) \bar{x}^{1-d}$ so that 
\begin{equation}
\label{ee 001}
\sup_{x\geq 0} f(x) = f(\bar{x}) \sim \frac{1}{bd} \bar{x}^{1-d} \exp(-b (b d)^{-1}) \sim \frac{1}{b d} (b d)^{-\frac{1-d}{d}} \exp(-1/d) = \frac{(d \exp(1))^{-1/d}}{b^{1/d}},
\end{equation}
whence the desired estimate \eqref{key bound}.

If $d<1$, we can verify that $p$ has two roots $x_0$ and $x_1$ on $\Rb_+$ with $x_0 < \bar{x}$ satisfying $x_0 \to 0$ as $b \downarrow 0$ and $x_1 \to +\infty $ as $b \downarrow 0$. Moreover we can check that $f$ admits a local minimum at $x_0$ and local maximum at $x_1$, so that the global maximum of $f$ on $[0,+\infty)$ is $\max(f(0),f(x_1)) = f(x_1)$ when $b \downarrow 0$. Similarly to the case $d\geq 1$, the roots $x_1$ satisfies $x_1^{d} \sim 1/bd$ and estimate \eqref{ee 001} is also valid for $\bar{x} = x_1$ whence \eqref{key bound}.

%\textbf{Acknowledgments.} The author would like to thank Rajinder Singh for inspiration.
\bibliographystyle{abbrv}
%\tiny
\bibliography{References_paper_TAO}

\begin{thebibliography}{10}

\bibitem{bailey1994fast}
D.~H. Bailey and P.~N. Swarztrauber.
\newblock A fast method for the numerical evaluation of continuous fourier and
  laplace transforms.
\newblock {\em SIAM Journal on Scientific Computing}, 15(5):1105--1110, 1994.

\bibitem{balakrishnan1985anomalous}
V.~Balakrishnan.
\newblock Anomalous diffusion in one dimension.
\newblock {\em Physica A: Statistical Mechanics and its applications},
  132(2-3):569--580, 1985.

\bibitem{chechkin2005fractional}
A.~V. Chechkin, R.~Gorenflo, and I.~M. Sokolov.
\newblock Fractional diffusion in inhomogeneous media.
\newblock {\em Journal of Physics A: Mathematical and General}, 38(42):L679,
  2005.

\bibitem{hao2019stability}
D.~N. H{\`a}o, J.~Liu, N.~Van~Duc, and N.~Van~Thang.
\newblock Stability results for backward time-fractional parabolic equations.
\newblock {\em Inverse Problems}, 35(12):125006, 2019.

\bibitem{kilbas2006theory}
A.~A. Kilbas, H.~M. Srivastava, and J.~J. Trujillo.
\newblock {\em Theory and applications of fractional differential equations},
  volume 204.
\newblock elsevier, 2006.

\bibitem{liu2010backward}
J.~Liu and M.~Yamamoto.
\newblock A backward problem for the time-fractional diffusion equation.
\newblock {\em Applicable Analysis}, 89(11):1769--1788, 2010.

\bibitem{louis1990mollifier}
A.~K. Louis and P.~Maass.
\newblock A mollifier method for linear operator equations of the first kind.
\newblock {\em Inverse problems}, 6(3):427, 1990.

\bibitem{mathe2008general}
P.~Math{\'e} and B.~Hofmann.
\newblock How general are general source conditions?
\newblock {\em Inverse Problems}, 24(1):015009, 2008.

\bibitem{metzler2000random}
R.~Metzler and J.~Klafter.
\newblock The random walk's guide to anomalous diffusion: a fractional dynamics
  approach.
\newblock {\em Physics reports}, 339(1):1--77, 2000.

\bibitem{podlubnv1999fractional}
I.~Podlubnv.
\newblock Fractional differential equations academic press.
\newblock {\em San Diego, Boston}, 1999.

\bibitem{schock1985approximate}
E.~Schock.
\newblock Approximate solution of ill-posed equations: arbitrarily slow
  convergence vs. superconvergence.
\newblock In {\em Constructive methods for the practical treatment of integral
  equations}, pages 234--243. Springer, 1985.

\bibitem{van2020mollification}
N.~Van~Duc, P.~Q. Muoi, and N.~Van~Thang.
\newblock A mollification method for backward time-fractional heat equation.
\newblock {\em Acta Mathematica Vietnamica}, 45(3):749--766, 2020.

\bibitem{wang2015optimal}
J.-G. Wang, T.~Wei, and Y.-B. Zhou.
\newblock Optimal error bound and simplified tikhonov regularization method for
  a backward problem for the time-fractional diffusion equation.
\newblock {\em Journal of computational and applied mathematics}, 279:277--292,
  2015.

\bibitem{wang2012data}
L.~Wang and J.~Liu.
\newblock Data regularization for a backward time-fractional diffusion problem.
\newblock {\em Computers \& Mathematics with Applications}, 64(11):3613--3626,
  2012.

\bibitem{wang2013total}
L.~Wang and J.~Liu.
\newblock Total variation regularization for a backward time-fractional
  diffusion problem.
\newblock {\em Inverse problems}, 29(11):115013, 2013.

\bibitem{xiong2012inverse}
X.~Xiong, Q.~Zhou, and Y.~Hon.
\newblock An inverse problem for fractional diffusion equation in 2-dimensional
  case: Stability analysis and regularization.
\newblock {\em Journal of Mathematical Analysis and Applications},
  393(1):185--199, 2012.

\bibitem{yang2014mollification}
F.~Yang, C.-L. Fu, and X.-X. Li.
\newblock A mollification regularization method for unknown source in
  time-fractional diffusion equation.
\newblock {\em International Journal of Computer Mathematics},
  91(7):1516--1534, 2014.

\bibitem{yang2013solving}
M.~Yang and J.~Liu.
\newblock Solving a final value fractional diffusion problem by boundary
  condition regularization.
\newblock {\em Applied Numerical Mathematics}, 66:45--58, 2013.

\bibitem{yang2015fourier}
M.~Yang and J.~Liu.
\newblock Fourier regularization for a final value time-fractional diffusion
  problem.
\newblock {\em Applicable Analysis}, 94(7):1508--1526, 2015.

\end{thebibliography}
\end{document}